\documentclass[10pt]{amsart}

\usepackage[T1]{fontenc}
\usepackage[nomath]{lmodern}
\usepackage{amsmath,amssymb,amsthm,mathtools}
\usepackage[letterpaper,margin=1in]{geometry}
\usepackage{microtype}
\DeclareMicrotypeSet[expansion]{textfonts}{encoding=T1,family={rm*,sf*}}
\UseMicrotypeSet[expansion]{textfonts}
\usepackage{aliascnt}
\usepackage[hidelinks]{hyperref}
\usepackage[nameinlink,capitalise]{cleveref}

\newtheorem{theorem}{Theorem}[section]
\newaliascnt{lemma}{theorem}
\newtheorem{lemma}[lemma]{Lemma}
\aliascntresetthe{lemma}
\newaliascnt{proposition}{theorem}
\newtheorem{proposition}[proposition]{Proposition}
\aliascntresetthe{proposition}
\newaliascnt{corollary}{theorem}
\newtheorem{corollary}[corollary]{Corollary}
\aliascntresetthe{corollary}
\theoremstyle{definition}
\newaliascnt{definition}{theorem}
\newtheorem{definition}[definition]{Definition}
\aliascntresetthe{definition}
\newaliascnt{example}{theorem}
\newtheorem{example}[example]{Example}
\aliascntresetthe{example}
\theoremstyle{remark}
\newaliascnt{remark}{theorem}
\newtheorem{remark}[remark]{Remark}
\aliascntresetthe{remark}

\crefname{theorem}{Theorem}{Theorems}
\Crefname{theorem}{Theorem}{Theorems}
\crefname{lemma}{Lemma}{Lemmas}
\Crefname{lemma}{Lemma}{Lemmas}
\crefname{proposition}{Proposition}{Propositions}
\Crefname{proposition}{Proposition}{Propositions}
\crefname{corollary}{Corollary}{Corollaries}
\Crefname{corollary}{Corollary}{Corollaries}
\crefname{definition}{Definition}{Definitions}
\Crefname{definition}{Definition}{Definitions}
\crefname{example}{Example}{Examples}
\Crefname{example}{Example}{Examples}
\crefname{remark}{Remark}{Remarks}
\Crefname{remark}{Remark}{Remarks}

\newcommand{\Open}{\operatorname{Open}}
\newcommand{\Left}{\operatorname{Left}}
\newcommand{\Inc}{\operatorname{Inc}}
\newcommand{\DES}{\operatorname{DES}}
\newcommand{\SYT}{\operatorname{SYT}}
\newcommand{\SSYT}{\operatorname{SSYT}}
\newcommand{\sgn}{\operatorname{sgn}}
\newcommand{\sink}{\operatorname{sink}}
\newcommand{\content}{\operatorname{content}}
\newcommand{\Acal}{\mathcal{A}}
\newcommand{\Mcal}{\mathcal{M}}
\newcommand{\Ccal}{\mathcal{C}}
\newcommand{\Tcal}{\mathcal{T}}
\newcommand{\CAcal}{\mathcal{CA}}
\newcommand{\frakA}{\mathfrak{A}}
\newcommand{\cyc}{\operatorname{cyc}}
\newcommand{\tr}{\operatorname{tr}}
\newcommand{\diag}{\operatorname{diag}}
\newcommand{\type}{\operatorname{type}}
\newcommand{\Stab}{\operatorname{Stab}}
\newcommand{\perm}{\operatorname{perm}}
\newcommand{\Wcal}{\mathcal{W}}
\newcommand{\Rcal}{\mathcal{R}}
\newcommand{\ind}{\mathbf{1}}
\newcommand{\SRHT}{\operatorname{SRHT}}
\newcommand{\rht}{\operatorname{ht}}
\newcommand{\Alt}{\operatorname{Alt}}
\newcommand{\Val}{\operatorname{Val}}
\newcommand{\rp}{\operatorname{rp}}
\newcommand{\pred}{\operatorname{pred}}
\numberwithin{equation}{section}
\hypersetup{pdftitle={Fixed-opening-set positivity for perfect matchings},pdfauthor={Myungkyu Chae},pdflang={en-US}}

\title[Fixed-opening-set positivity]{Fixed-opening-set positivity for perfect matchings}
\author{Myungkyu Chae}
\address{Department of Mathematics, Washington University in St. Louis, St. Louis, Missouri 63130, USA}
\email{c.myungkyu@wustl.edu}
\date{}

\begin{document}
\begin{abstract}
For a perfect matching $M$ on $[2n]$, let $P_M$ be the interval order on its edges and let
$\Open(M)$ be the set of smaller endpoints of its edges. For an opening set $O$ that occurs for
some matching, put $S_O=\sum_{\Open(M)=O}X_{\Inc(P_M)}$ and $H_O=\omega S_O$.
We give explicit Schur and power-sum expansions of $H_O$ with nonnegative integer coefficients,
and prove that $H_O$ is $h$-positive (equivalently, $S_O$ is $e$-positive).
The Schur coefficients count standard tableaux of staggered skew shapes with prescribed entries in
the left cells of distinguished dominoes. The Schur expansion follows from a Jacobi--Trudi sign-reversing
involution on arrays with prescribed opener positions.
For every partition $\mu$ of $n$, a coefficient-transfer formula expresses $[h_\mu]H_O$ as a
positive scalar multiple of an elementary coefficient of the chromatic symmetric function of a
clique blow-up of a path. This graph is a unit interval graph, so Hikita's theorem gives
nonnegativity. The power-sum coefficients admit a cycle-type autocorrelation formula in terms of
the restricted permutations from the Ferrers threshold model and an interpretation by alternating
permutations with prescribed valley values and record partitions; we also obtain explicit threshold formulas.
Further results include descent-set refinements, formulas for selected $h$-coefficients, and a
rim-hook cancellation for the coefficient of $h_n$. These results refine the global matching identity of Amdeberhan,
Shareshian, and Stanley.
\end{abstract}
\maketitle
\pagestyle{plain}

\section{Introduction}

A perfect matching $M$ on $[2n]$ defines an interval order $P_M$ on its edges: for $a<b$ and $c<d$, set
\[
 \{a,b\}<_{P_M}\{c,d\}\quad\Longleftrightarrow\quad b<c.
\]
Write $\Inc(P_M)$ for the incomparability graph of $P_M$: its vertices are the edges of $M$,
and two distinct vertices are adjacent exactly when they are incomparable in $P_M$.
Symmetric-function conventions, including $\omega$, are given in \cref{subsec:symmetric-functions};
the matching and incomparability-graph notation is defined in \cref{subsec:matchings}.
Amdeberhan, Shareshian, and Stanley \cite[Theorem~7.1]{ASS} proved the identity
\begin{equation}\label{eq:global-sprout}
 (2n)!A_n=\sum_{M\in\Mcal(2n)}\omega X_{\Inc(P_M)},
\end{equation}
where $\Mcal(2n)$ is the set of perfect matchings on $[2n]$, $X_G$ is Stanley's chromatic symmetric
function \cite{StanleyXG}, and
\[
 \sum_{n\ge0}A_n(x)t^n=\prod_{i\ge1}\sec\!\bigl(\sqrt{x_it}\bigr).
\]
They showed that \eqref{eq:global-sprout} is $h$-positive and described its Schur coefficients using
standard Young tableaux of staggered skew shapes \cite[Theorems~5.1 and~6.2]{ASS}.

We refine these expansions by prescribing the set of smaller endpoints of the matching edges, called
\emph{openers}. For $O\subseteq[2n]$, let $\Mcal(O)$ be the set of perfect matchings with opening set
$O$, and call $O$ valid if $\Mcal(O)\ne\varnothing$ (see \cref{subsec:matchings}). Put
\[
 S_O=\sum_{M\in\Mcal(O)}X_{\Inc(P_M)},\qquad H_O=\omega S_O.
\]
Our main result is
\begin{equation}\label{eq:intro-h-positive}
 H_O=\sum_{\mu\vdash n}c^O_\mu h_\mu,
 \qquad c^O_\mu\in\mathbb Z_{\ge0},
\end{equation}
for every valid $O$. Equivalently, $S_O$ is $e$-positive. Individual summands need not even be
Schur-positive: the matching $18\mid23\mid45\mid67$ has incomparability graph $K_{1,3}$, the four-vertex star
with center $\{1,8\}$ and leaves $\{2,3\},\{4,5\},\{6,7\}$; $X_{K_{1,3}}$ is not Schur-positive
\cite[p.~186]{StanleyXG}.
Indeed, $\{1,8\}$ is incomparable with the other three edges, whereas
$\{2,3\}<_{P_M}\{4,5\}<_{P_M}\{6,7\}$.

The proof uses an auxiliary graph determined by $O$. With the threshold notation of
\cref{subsec:thresholds}, write $O=\{o_1<\cdots<o_n\}$ and set
\[
 r_i=o_i-i+1,\qquad b_t=\#\{i:r_i=t\},\qquad D_O=\prod_{t=1}^n b_t!.
\]
Form a word of length $2n$ by writing $O$ at positions belonging to the opening set and $C$ at
the remaining positions. Reading $O$ as an up-step and $C$ as a down-step gives a Dyck path: it
starts at height $0$, never falls below height $0$, and ends at height $0$. Let $\eta_j(O)$ be its
height immediately before the $j$th down-step, equivalently, the $j$th closer. Define
\[
 \gamma_{2j-1}=n-\eta_j(O)\quad(1\le j\le n),\qquad
 \gamma_{2j}=\eta_j(O)-1\quad(1\le j<n).
\]
Let $G_O$ be the clique blow-up of the path $1-2-\cdots-(2n-1)$: replace each vertex $s$ by a
clique of size $\gamma_s$, and join every vertex in one clique to every vertex in the next clique.
There are no other edges between cliques. A zero size leaves that position empty, with no edge
bridging the gap. The formal construction is given in \cref{subsec:auxiliary-graph}.
For $\mu\vdash n$, padded with zeros to length $n$, put
\[
 \nu_\mu=(n-\mu_n,n-\mu_{n-1},\ldots,n-\mu_1),
 \qquad \gamma!=\prod_{s=1}^{2n-1}\gamma_s!,
\]
omitting zero parts from $\nu_\mu$. The coefficient-transfer theorem (\cref{thm:coefficient-transfer}) gives
\begin{equation}\label{eq:intro-transfer}
 [h_\mu]H_O=\frac{D_O}{\gamma!}[e_{\nu_\mu}]X_{G_O}
 \qquad(\mu\vdash n).
\end{equation}
The graph $G_O$ has $n(n-1)$ vertices and is a unit interval graph (\cref{subsec:unit-interval}). Hikita's theorem at $q=1$
\cite{Hikita2025SS}, which resolves the Stanley--Stembridge conjecture \cite{StanleyStembridge},
therefore gives \eqref{eq:intro-h-positive}. The proof of the transfer formula is carried out in \cref{subsec:repeated-thresholds,subsec:inverse-path,subsec:rectangle-complement}:
repeated thresholds are first grouped in a universal matrix, the inverse matrix has path support, and a
rectangle-complement identity then produces the coefficient of $X_{G_O}$.

The Schur coefficients in $H_O$ have a direct tableau interpretation. For $\lambda\vdash n$, let
$\rho(\lambda)$ have row lengths $2\lambda'_1,2\lambda'_2,\ldots$, with each row below the first shifted one
cell to the left of the row above it. Divide each row into consecutive horizontal dominoes, and let $\Left(A)$ be the
set of entries in the left cells of these dominoes. Write $\SYT(\rho(\lambda))$ for the standard Young
tableaux of this shape, as in \cref{subsec:tableaux}. With
\[
 g_O^\lambda=\#\{A\in\SYT(\rho(\lambda)):\Left(A)=O\},
\]
we prove
\begin{equation}\label{eq:intro-schur}
 S_O=\sum_{\lambda\vdash n}g_O^\lambda s_{\lambda'},
 \qquad H_O=\sum_{\lambda\vdash n}g_O^\lambda s_\lambda.
\end{equation}

The power-sum coefficients of $H_O$ have two positive combinatorial descriptions. The threshold
bijection of \cref{prop:threshold} identifies
$P_r=\{\pi\in\mathfrak S_n:\pi(i)\ge r_i\text{ for all }i\}$ with $\Mcal(O)$. With
$\Alt_{2n}(O)$ denoting the alternating permutations with valley-value set $O$, we prove
\begin{equation}\label{eq:intro-power-sum}
 H_O=\sum_{\pi,\rho\in P_r}p_{\cyc(\rho^{-1}\pi)}
 =\sum_{w\in\Alt_{2n}(O)}p_{\rp(\widehat w)}.
\end{equation}
Here $\cyc(\tau)$ is the partition of the cycle lengths of $\tau$, including fixed points.
The word $\widehat w$ consists of the odd-position entries of $w$, and $\rp$ denotes its record
partition; these statistics and the alternating convention are defined in
\cref{subsec:fixed-valley-power-sum}.
The alternating-permutation formula refines the record-partition interpretation of
\cite[Theorem~3.2]{ASS} by fixing the valley-value set; summing over all valid $O$
recovers the power-sum expansion of $(2n)!A_n$ (\cref{cor:global-record-refinement}).
In particular, $H_O$ is $p$-positive. This also follows from \eqref{eq:intro-h-positive}, since
$h_k=\sum_{\lambda\vdash k}p_\lambda/z_\lambda$. We also give explicit formulas for these
coefficients in terms of the threshold sequence $r$.

Fixed opener and closer sets occur in the study of crossings and nestings
\cite{ChenDengDuStanleyYan}; related Schur-positive matching statistics appear in
\cite{AdinRoichman,Marmor}. Our other results include descent-set identities, formulas for
selected complete homogeneous coefficients, and a rim-hook involution for the coefficient of $h_n$. The appendix
records a common Ferrers-board model for the stable partitions of all matchings in $\Mcal(O)$.

\section{Preliminaries and notation}\label{sec:prelim}

Write $[d]=\{1,\ldots,d\}$ for $d\ge1$, and $[0]=\varnothing$.
Throughout, $n\ge1$. All graphs are finite and simple unless a directed graph or multigraph is explicitly
introduced. A partition of $n$ is a weakly decreasing sequence of positive integers with sum $n$,
denoted by $\lambda\vdash n$. Its conjugate is $\lambda'$, where
$\lambda'_j=\#\{i:\lambda_i\ge j\}$, its size is $|\lambda|=n$, and its number of positive parts
is $\ell(\lambda)$. A composition of $n$ is an ordered sequence of positive integers with sum $n$;
a weak composition allows zero parts. Products over empty index sets are $1$.
The symmetric group on $[d]$ is denoted by $\mathfrak S_d$. Permutations are composed from right to left;
one-line notation $\pi(1)\cdots\pi(d)$ lists their values, and $\sgn(\pi)$ denotes their sign.
For a permutation $\tau$, $\cyc(\tau)$ is the partition of its cycle lengths, including fixed points.

\subsection{Symmetric functions and the Hall inner product}\label{subsec:symmetric-functions}

Let $\Lambda_{\mathbb Z}=\mathbb Z[e_1,e_2,\ldots]$ be the ring of symmetric functions with integer
coefficients and let $\Lambda_{\mathbb Q}=\Lambda_{\mathbb Z}\otimes_{\mathbb Z}\mathbb Q$.
We use the standard monomial, complete homogeneous, elementary, power-sum, and Schur bases
$m_\lambda$, $h_\lambda$, $e_\lambda$, $p_\lambda$, and $s_\lambda$, respectively. Here
\[
 h_k=\sum_{i_1\le\cdots\le i_k}x_{i_1}\cdots x_{i_k},\qquad
 e_k=\sum_{i_1<\cdots<i_k}x_{i_1}\cdots x_{i_k},\qquad p_k=\sum_i x_i^k,
\]
where the sums range over positive indices, and $h_\lambda=\prod_i h_{\lambda_i}$, $e_\lambda=\prod_i e_{\lambda_i}$,
and $p_\lambda=\prod_i p_{\lambda_i}$. The function $m_\lambda$ is the sum of the distinct monomials
whose positive exponents are the parts of $\lambda$; the tableau description of $s_\lambda$ is recalled
in \cref{subsec:tableaux}. For any basis $\{b_\lambda\}$, $[b_\lambda]f$ denotes the corresponding
coefficient of $f$, and $b$-positive means that all these coefficients are nonnegative. Brackets such
as $[x_1^{a_1}\cdots x_d^{a_d}]F$ similarly extract an ordinary monomial coefficient.
The $m$-, $h$-, $e$-, and $s$-bases are integral bases of $\Lambda_{\mathbb Z}$, and the $p$-basis
is a basis over $\mathbb Q$. Set $h_0=e_0=1$ and $h_k=e_k=0$ for $k<0$. The standard involution
$\omega$ is the algebra automorphism with $\omega^2=\mathrm{id}$ satisfying
$\omega(e_k)=h_k$, $\omega(s_\lambda)=s_{\lambda'}$, and $\omega(p_k)=(-1)^{k-1}p_k$. The
Hall inner product is normalized so that the monomial and complete homogeneous bases are dual and
the Schur basis is orthonormal:
\begin{equation}\label{eq:Hall-duality}
\langle m_\lambda,h_\mu\rangle=\delta_{\lambda\mu},
\qquad
\langle s_\lambda,s_\mu\rangle=\delta_{\lambda\mu}.
\end{equation}
Here $\delta_{\lambda\mu}$ is the Kronecker delta: it is $1$ when $\lambda=\mu$ and $0$ otherwise.
Thus orthonormality means that distinct Schur functions are orthogonal and each has squared norm $1$.
Also,
\begin{equation}\label{eq:Hall-power}
 \langle p_\lambda,p_\mu\rangle=\delta_{\lambda\mu}z_\lambda,
 \qquad z_\lambda=\prod_{j\ge1}j^{m_j(\lambda)}m_j(\lambda)!,
\end{equation}
where $m_j(\lambda)$ denotes the multiplicity of the part $j$. In particular,
$[h_\mu]f=\langle f,m_\mu\rangle$. For an integer partition $\lambda$, define the augmented monomial function by
\[
 \widetilde m_\lambda=\left(\prod_{j\ge1}m_j(\lambda)!\right)m_\lambda.
\]
These conventions follow \cite[Chapter~7]{StanleyEC2}. For example,
\[
 \omega(e_{(2,2)})=h_2^2,\qquad
 \langle s_{(2,2)},s_{(2,2)}\rangle=1,\qquad
 \langle s_{(2,2)},s_{(3,1)}\rangle=0.
\]

\subsection{Perfect matchings and interval orders}\label{subsec:matchings}

Let $\Mcal(2n)$ denote the set of all perfect matchings on $[2n]$.
A perfect matching is a partition of $[2n]$ into $n$ two-element blocks, called edges.
A perfect matching $M$ is written as a set of edges $\{a,b\}$ with $a<b$. We call $a$ its opener
and $b$ its closer, and abbreviate $\{a,b\}$ by $ab$ in displayed matching lists. Define its opening set by
\begin{equation}\label{eq:opening-set}
\Open(M):=\{a:\{a,b\}\in M,\ a<b\},
\end{equation}
and let
\begin{equation}\label{eq:McalO}
\Mcal(O):=\{M:\Open(M)=O\}.
\end{equation}
A set $O\subseteq[2n]$ is \emph{valid} if $|O|=n$ and $\Mcal(O)\ne\varnothing$.
For an $n$-element set $O\subseteq[2n]$, write $C=[2n]\setminus O$. Then $O$ is valid if and only if
\begin{equation}\label{eq:valid-prefix}
 |O\cap[k]|\ge |C\cap[k]|\qquad(0\le k\le 2n).
\end{equation}
Indeed, every closer must be paired with an earlier opener. Conversely, under
\eqref{eq:valid-prefix}, a left-to-right scan can pair each closer with an unmatched earlier opener,
and all vertices are paired by the end. Equivalently, if $O=\{o_1<\cdots<o_n\}$, then $o_i\le 2i-1$
for all $i$. We henceforth fix a valid $O$.
The word $w(O)\in\{O,C\}^{2n}$ has letter $O$ at positions in the opening set and letter $C$ at the
remaining positions. By \eqref{eq:valid-prefix}, it is a Dyck word: reading $O$ as an up-step
and $C$ as a down-step gives a path that stays at nonnegative height and ends at height $0$.
The height of a prefix is its number of $O$ letters minus its number of $C$ letters.

For $M\in\Mcal(O)$, define the interval order $P_M$ on the edges of $M$ by
\begin{equation}\label{eq:interval-order}
\{a,b\}<_{P_M}\{c,d\}
\quad\Longleftrightarrow\quad
b<c.
\end{equation}
We write $\le_{P_M}$ for the reflexive closure of $<_{P_M}$, and $\Inc(P_M)$ for the incomparability
graph of $P_M$: its vertices are the matching edges, and two distinct vertices are adjacent
exactly when neither precedes the other in $P_M$. A chain consists of pairwise comparable
elements. A stable (or independent) set of a graph contains no adjacent pair of vertices, and
a clique consists of pairwise adjacent vertices. This is the interval-order construction of
\cite[Section~7]{ASS}.

For a finite graph $G$, Stanley's chromatic symmetric function is
\[
 X_G(x)=\sum_{\kappa\text{ proper}}\prod_{v\in V(G)}x_{\kappa(v)},
\]
where a proper coloring $\kappa:V(G)\to\mathbb Z_{>0}$ assigns different colors to adjacent vertices
\cite{StanleyXG}. We set $X_\varnothing=1$.  For fixed $O$, set
\begin{equation}\label{eq:SO-def}
S_O:=\sum_{M\in\Mcal(O)}X_{\Inc(P_M)},
\qquad
H_O:=\omega S_O.
\end{equation}
If
\[
S_O=\sum_{\lambda\vdash n}c^O_\lambda e_\lambda,
\]
then, since $\omega(e_\lambda)=h_\lambda$,
\begin{equation}\label{eq:h-e-duality}
H_O=\sum_{\lambda\vdash n}c^O_\lambda h_\lambda,
\qquad
[h_\lambda]H_O=[e_\lambda]S_O.
\end{equation}

A \emph{sequencing} of $M$ is an ordering
\[
\sigma=\sigma_1\sigma_2\cdots\sigma_n
\]
of the edges of $M$.  Following \cite[Theorem~7.1]{ASS}, define
\begin{equation}\label{eq:des-sequencing}
\DES(\sigma)
:=
\{i\in[n-1]:\sigma_i\not\le_{P_M}\sigma_{i+1}\}.
\end{equation}
Since consecutive terms of a sequencing are distinct, $\not\le_{P_M}$ in \eqref{eq:des-sequencing}
is equivalent to $\not<_{P_M}$, agreeing with Chow's convention \cite[Corollary~2]{Chow}. For
$D,T\subseteq[n-1]$, set
\begin{align}
N_O(D)
&:=
\#\{(M,\sigma):M\in\Mcal(O),\ \DES(\sigma)=D\},
\label{eq:NO}\\
R_O(T)
&:=
\sum_{D\subseteq T}N_O(D)
=
\#\{(M,\sigma):M\in\Mcal(O),\ \DES(\sigma)\subseteq T\}.
\label{eq:RO}
\end{align}

If
\[
T=\{t_1<\cdots<t_k\}\subseteq[n-1],
\qquad
t_0=0,
\qquad
t_{k+1}=n,
\]
define the associated composition
\begin{equation}\label{eq:alphaT}
\alpha(T):=
(t_1-t_0,t_2-t_1,\ldots,t_{k+1}-t_k).
\end{equation}

\begin{example}[Running example]\label{ex:running-matching}
Let $n=4$, $O=\{1,2,4,6\}$, and $C=\{3,5,7,8\}$. Then $w(O)=OOCOCOCC$, and
$M=13\mid25\mid47\mid68$ belongs to $\Mcal(O)$. The sequencing
$\sigma=(13,47,25,68)$ has $\DES(\sigma)=\{2\}$, since
$13<_{P_M}47$ and $25<_{P_M}68$, while $47\not<_{P_M}25$. Thus
$\alpha(\{2\})=(2,2)$. We use this opening set throughout the examples below.
\end{example}

\subsection{Quasisymmetric functions}\label{subsec:quasisymmetric}

For $S\subseteq[n-1]$, define the fundamental quasisymmetric function
\begin{equation}\label{eq:fundamental-qsy}
Q_{S,n}
:=
\sum_{\substack{1\le i_1\le\cdots\le i_n\\
j\in S\Rightarrow i_j<i_{j+1}}}
x_{i_1}\cdots x_{i_n},
\end{equation}
and the monomial quasisymmetric function in subset notation by
\begin{equation}\label{eq:monomial-qsy}
M_{S,n}
:=
\sum_{\substack{1\le i_1\le\cdots\le i_n\\
i_j<i_{j+1}\iff j\in S}}
x_{i_1}\cdots x_{i_n}.
\end{equation}
Thus $M_{S,n}$ records the exact strict-rise set, whereas $Q_{S,n}$ only requires strict rises at
the positions in $S$.  Consequently,
\begin{equation}\label{eq:Q-to-M}
Q_{D,n}
=
\sum_{D\subseteq T\subseteq[n-1]}M_{T,n}.
\end{equation}
This is the subset-indexed transition between the fundamental and monomial quasisymmetric bases;
compare \cite[Section~7.19]{StanleyEC2} and \cite{Chow}. For the set $T=\{2\}$ from
\cref{ex:running-matching}, $M_{\{2\},4}=\sum_{i<j}x_i^2x_j^2=m_{(2,2)}$, and
\[
 Q_{\{2\},4}=M_{\{2\},4}+M_{\{1,2\},4}+M_{\{2,3\},4}+M_{\{1,2,3\},4}.
\]

\subsection{Tableaux}\label{subsec:tableaux}

The Young diagram of $\lambda$ has $\lambda_i$ cells in row $i$; a skew shape $\beta/\delta$
is obtained by removing the diagram of $\delta$ from that of $\beta$. For a partition
$\lambda\vdash n$, let $\rho(\lambda)$ be the staggered skew shape of
\cite[Section~6]{ASS}.  Its row lengths are
\[
2\lambda'_1,2\lambda'_2,\ldots,
\]
and each lower row begins one cell to the left of the row immediately above it. More explicitly,
with $d=\ell(\lambda')$,
\begin{equation}\label{eq:rho-shape}
 \rho(\lambda)=\beta/\delta,
 \qquad \beta_i=2\lambda'_i+d-i,\qquad \delta_i=d-i
 \quad(1\le i\le d).
\end{equation}
We use English notation, with rows numbered from top to bottom. A standard Young tableau of a
(possibly skew) shape $\theta$ with $|\theta|$ cells fills them with $1,\ldots,|\theta|$, each used once,
increasing along rows from left to right and columns from top to bottom; the set of all such tableaux
is denoted by $\SYT(\theta)$. Each row of $\rho(\lambda)$
is divided into consecutive pairs starting at its leftmost cell; these pairs are the distinguished
horizontal dominoes. If $A\in\SYT(\rho(\lambda))$, let $\Left(A)$ be the set of entries in the left
cells of these dominoes. Pairing the two entries of every distinguished domino produces a perfect matching, so
$\Left(A)$ is always a valid opening set.  Set
\begin{equation}\label{eq:gOlambda}
g_O^\lambda
:=
\#\{A\in\SYT(\rho(\lambda)):\Left(A)=O\}.
\end{equation}

Write $\SSYT(\mu)$ for the semistandard Young tableaux of shape $\mu$, filled with positive integers,
whose rows are weakly increasing and columns strictly increasing. Let $\SSYT(\mu,\alpha)$ be those of content $\alpha$:
symbol $i$ occurs $\alpha_i$ times. The Schur function is
$s_\mu(x)=\sum_{U\in\SSYT(\mu)}x^{\content(U)}$, where
$x^\alpha=\prod_i x_i^{\alpha_i}$. For a partition $\mu\vdash n$ and a composition $\alpha$ of $n$,
let
\begin{equation}\label{eq:Kostka-composition}
K_{\mu,\alpha}
:=
\#\{U\in\SSYT(\mu):\content(U)=\alpha\}.
\end{equation}
If $Q\in\SYT(\mu)$, define
\[
\DES(Q)
:=
\{i\in[n-1]: i+1\text{ lies in a lower row than }i\},
\]
and set
\begin{equation}\label{eq:fDmu}
f_D^\mu
:=
\#\{Q\in\SYT(\mu):\DES(Q)=D\}.
\end{equation}

\begin{example}\label{ex:running-tableau}
For $O=\{1,2,4,6\}$ and $\lambda=(2,2)$, we have $\rho(\lambda)=(5,4)/(1)$, and
\[
 A=\begin{array}{ccccc}
  &\mathbf{1}&3&\mathbf{4}&7\\
  \mathbf{2}&5&\mathbf{6}&8&
 \end{array}
 \qquad\text{satisfies}\qquad \Left(A)=O.
\]
The bold entries are the distinguished left entries. Pairing consecutive cells in each row
recovers the matching in \cref{ex:running-matching}. For the ordinary shape $(2,2)$, the unique
semistandard tableau of content $(2,2)$ has rows $(1,1)$ and $(2,2)$, so $K_{(2,2),(2,2)}=1$.
The standard tableau with rows $(1,2)$ and $(3,4)$ has descent set $\{2\}$.
\end{example}

\subsection{Thresholds and pre-closure heights}\label{subsec:thresholds}

Write
\[
 O=\{o_1<\cdots<o_n\},\qquad C=\{c_1<\cdots<c_n\},
\]
and define
\begin{equation}\label{eq:threshold-r}
 r_i=1+\#\{j:c_j<o_i\}=o_i-i+1.
\end{equation}
Validity gives $o_1=1$, $o_i\le2i-1$, and consequently
\begin{equation}\label{eq:threshold-bounds}
 1=r_1\le r_2\le\cdots\le r_n,\qquad r_i\le i.
\end{equation}
The monotonicity follows from $r_{i+1}-r_i=o_{i+1}-o_i-1\ge0$. Conversely, every integer sequence
satisfying \eqref{eq:threshold-bounds} comes from a unique valid opening set: $o_i=r_i+i-1$ is
strictly increasing and satisfies $o_i\le2i-1$. We call such an $r$ a valid threshold sequence. Thus $r_i$ is the first closer rank allowed
for the opener $o_i$. Put
\begin{equation}\label{eq:Pr}
 P_r=\{\pi\in\mathfrak S_n:\pi(i)\ge r_i\text{ for every }i\in[n]\}.
\end{equation}

\begin{proposition}[Ferrers threshold encoding]\label{prop:threshold}
The map
\begin{equation}\label{eq:matching-from-pi}
 \pi\longmapsto M_\pi:=\{\{o_i,c_{\pi(i)}\}:i\in[n]\}
\end{equation}
is a bijection from $P_r$ to $\Mcal(O)$.
\end{proposition}
\begin{proof}
Precisely $c_1,\ldots,c_{r_i-1}$ precede $o_i$, so $o_i<c_j$ if and only if $j\ge r_i$. A
permutation pairs every opener with a distinct closer, and this threshold condition is exactly what
makes each specified opener the smaller endpoint. The inverse reads the closer rank paired with each
opener.
\end{proof}

Define
\begin{equation}\label{eq:threshold-data}
 b_t=\#\{i:r_i=t\},\qquad
 B_j=\sum_{t=1}^j b_t,\qquad
 D_O=\prod_{t=1}^n b_t!.
\end{equation}
Thus $b_t$ counts the openers with threshold $t$, and $B_j$ counts the openers whose threshold
is at most $j$. Since $r_i\le i$, we have $B_j\ge j$. Also $r_i\le j$ if and only if $o_i<c_j$, and hence the
pre-closure height is
\begin{equation}\label{eq:eta}
 \eta_j(O)=|O\cap[c_j-1]|-|C\cap[c_j-1]|
          =B_j-j+1.
\end{equation}
Thus
\begin{equation}\label{eq:eta-basic}
 1\le\eta_j(O)\le n,\qquad \eta_n(O)=1.
\end{equation}
We sometimes abbreviate $\eta_j(O)$ to $\eta_j$, and write
$\eta(O)=(\eta_1(O),\ldots,\eta_n(O))$.

\begin{lemma}\label{lem:matchingheight}
For every valid opening set $O$,
\begin{equation}\label{eq:matchingheight}
 |\Mcal(O)|=|P_r|=\prod_{j=1}^n\eta_j(O)
                     =\prod_{j=1}^{n-1}\eta_j(O).
\end{equation}
\end{lemma}
\begin{proof}
Process the closers in increasing order. Before $c_j$ is matched, $B_j$ openers have appeared and
$j-1$ of them have been used. Thus there are exactly $B_j-j+1=\eta_j(O)$ choices, for every partial
matching obtained so far. The multiplication principle and $\eta_n(O)=1$ give the formula.
\end{proof}

\begin{example}\label{ex:running-thresholds}
For $O=\{1,2,4,6\}$, the threshold and height data are
\[
 r=(1,1,2,3),\quad b=(2,1,1,0),\quad B=(2,3,4,4),\quad \eta(O)=(2,2,2,1).
\]
For instance, three openers precede $c_2=5$, and one closer has already occurred, so
$B_2=3$ and $\eta_2=3-1=2$.
Consequently $D_O=2!1!1!0!=2$ and $|\Mcal(O)|=2\cdot2\cdot2\cdot1=8$.
The identity permutation encodes $M=13\mid25\mid47\mid68$.
\end{example}

\section{Ferrers encoding and squarefree extraction}\label{sec:ferrers}

\subsection{Relabeling by closer rank}

We first regroup the matchings according to the threshold seen at each closer rank. For
$\pi\in P_r$, label the edges of $M_\pi$ by their closer ranks: vertex $v$ denotes the edge
$\{o_{\pi^{-1}(v)},c_v\}$. Set
\begin{equation}\label{eq:threshold-word}
 a_v^\pi=r_{\pi^{-1}(v)}.
\end{equation}
This is a rearrangement of $r$ with $a_v^\pi\le v$. Let
\[
 \Wcal(r)=\{a=(a_1,\ldots,a_n):a\text{ is a rearrangement of }r,
                            \ a_v\le v\text{ for all }v\}.
\]
For $a\in\Wcal(r)$, define $G_a$ on $[n]$ by
\begin{equation}\label{eq:Ga}
 \{u,v\}\in E(G_a),\quad u<v,
 \qquad\Longleftrightarrow\qquad u\ge a_v.
\end{equation}

\begin{proposition}[Graph identification and fiber multiplicity]\label{prop:regrouping}
Under the closer-rank labeling, $\Inc(P_{M_\pi})=G_{a^\pi}$. Each $a\in\Wcal(r)$ is attained by
exactly $D_O$ permutations in $P_r$. Consequently, writing
\[
 F_r=\sum_{a\in\Wcal(r)}X_{G_a},
\]
we have
\begin{equation}\label{eq:global-regrouping}
 S_O=D_OF_r=D_O\sum_{a\in\Wcal(r)}X_{G_a}.
\end{equation}
\end{proposition}
\begin{proof}
Let $u<v$. The edge with closer $c_v$ cannot lie strictly to the left of the edge with closer $c_u$.
The only possible comparability is therefore
\[
 \{o_{\pi^{-1}(u)},c_u\}<_{P_{M_\pi}}
 \{o_{\pi^{-1}(v)},c_v\}
 \quad\Longleftrightarrow\quad
 c_u<o_{\pi^{-1}(v)}
 \quad\Longleftrightarrow\quad u<a_v^\pi.
\]
The graph identification follows.

For fixed $a\in\Wcal(r)$, let $I_t=\{i:r_i=t\}$ and $V_t=\{v:a_v=t\}$. A permutation yields $a$
exactly when it restricts to a bijection $I_t\to V_t$ for every $t$. Each such choice is admissible,
because $r_i=t=a_{\pi(i)}\le\pi(i)$. The independent choices give $\prod_t b_t!=D_O$ permutations.
In particular every $a$ is attained, and regrouping the sum for $S_O$ proves
\eqref{eq:global-regrouping}.
\end{proof}

\begin{example}\label{ex:threshold-words}
For the data in \cref{ex:running-thresholds}, the admissible threshold words and their graph
edge sets are
\[
\begin{array}{c|l}
a&E(G_a)\\ \hline
(1,1,2,3)&\{\{1,2\},\{2,3\},\{3,4\}\}\\
(1,1,3,2)&\{\{1,2\},\{2,4\},\{3,4\}\}\\
(1,2,1,3)&\{\{1,3\},\{2,3\},\{3,4\}\}\\
(1,2,3,1)&\{\{1,4\},\{2,4\},\{3,4\}\}.
\end{array}
\]
The first two graphs are paths on four vertices and the last two are stars $K_{1,3}$.
The permutations $\pi=1234$ and $\pi'=2134$, written in one-line notation, both give the word
$a=(1,1,2,3)$. Their matchings are
\[
 M_\pi=13\mid25\mid47\mid68,\qquad
 M_{\pi'}=15\mid23\mid47\mid68.
\]
Under closer-rank labeling, both incomparability graphs are the path $1-2-3-4$.
Each threshold word has exactly $D_O=2$ preimages, obtained by interchanging the two openers with
threshold $1$. Consequently,
\[
 S_O=4X_{P_4}+4X_{K_{1,3}},
\]
where $P_4$ is the four-vertex path.
\end{example}

\subsection{Matrix evaluation and chain minors}
For an $n\times n$ matrix $A$ over a commutative $\mathbb Q$-algebra,
define the algebra homomorphism $f\mapsto f[A]$ by
\begin{equation}\label{eq:fos-matrix-evaluation}
\sum_{k\ge0}e_k[A]t^k=\det(I+tA).
\end{equation}
Thus $e_k[A]=0$ for $k>n$, and Newton's identities give
\begin{equation}\label{eq:matrix-trace}
p_k[A]=\tr(A^k).
\end{equation}
Here $I$ is the identity matrix and $\tr$ denotes the trace. The value $f[A]$ lies in the
coefficient algebra of $A$. Let
\begin{equation}\label{eq:fos-Ferrers-matrix}
\begin{aligned}
R_{ij}&:=\ind[j\ge r_i],&\qquad U&:=\diag(u_1,\ldots,u_n),\\
T&:=URVR^{\mathsf T},& V&:=\diag(v_1,\ldots,v_n),
\end{aligned}
\end{equation}
where all $u_i,v_i$ are independent commuting variables, $\ind[P]$ is $1$ if $P$ holds and $0$
otherwise, $\diag$ denotes a diagonal matrix, and ${}^{\mathsf T}$ denotes transpose. Thus $R$ is
the $0$--$1$ incidence matrix of the allowed opener--closer pairs $o_i<c_j$. In the running example,
the four rows begin with their first $1$ in columns $1,1,2,3$, respectively. Write
\[
\mathcal E F:=[u_1\cdots u_n v_1\cdots v_n]F.
\]
Thus $\mathcal E$ extracts the coefficient of the squarefree monomial in which every $u_i$ and every $v_i$ occurs once.
For index sets $I,J$, $R[I,J]$ denotes the corresponding submatrix.
All row and column indices of a submatrix are listed in increasing order.

\begin{lemma}[Ferrers minors and chains]\label{lem:chain-minors}
Let $I=\{i_1<\cdots<i_k\}$ and $J=\{j_1<\cdots<j_k\}$ be subsets of
$[n]$.  Then
\begin{equation}\label{eq:chain-minors}
\det R[I,J]^2=
\begin{cases}
1,&o_{i_1}<c_{j_1}<o_{i_2}<c_{j_2}<\cdots<o_{i_k}<c_{j_k},\\
0,&\text{otherwise}.
\end{cases}
\end{equation}
In the first case, the corresponding unique chain of matching edges is
$\{o_{i_1},c_{j_1}\}<\cdots<\{o_{i_k},c_{j_k}\}$.
The empty minor has determinant $1$.
\end{lemma}
\begin{proof}
The rows of $R[I,J]$ are weakly nested terminal strings of ones.
A nonzero determinant requires that these rows be distinct and nonzero.
There are exactly $k$ nonzero terminal strings of length $k$, so this
happens precisely when row $a$ has its first one in column $a$ for every
$a\in[k]$.  The resulting matrix is upper triangular with diagonal
entries $1$.  Equivalently,
$r_{i_a}\le j_a$ for all $a$ and $j_{a-1}<r_{i_a}$ for $a\ge2$.
By the definition of the thresholds, these inequalities say exactly
$o_{i_a}<c_{j_a}$ and $c_{j_{a-1}}<o_{i_a}$.  They prove
\eqref{eq:chain-minors}.  A chain with these opener and closer sets
must pair them in their respective increasing orders, proving uniqueness.
\end{proof}

Principal-minor expansion followed by Cauchy--Binet gives
\begin{equation}\label{eq:elementary-chain-evaluation}
e_k[T]=
\sum_{\substack{I,J\subseteq[n]\\|I|=|J|=k}}
\det R[I,J]^2\,u_Iv_J,
\end{equation}
where $u_I:=\prod_{i\in I}u_i$ and $v_J:=\prod_{j\in J}v_j$.
Thus $e_k[T]$ records exactly the possible chains of $k$ matching edges,
together with their opener and closer sets.

\begin{proposition}[Coloring kernel and squarefree extraction]
\label{prop:fos-squarefree}
In the degreewise completion of the ring of symmetric functions
(with identities interpreted one homogeneous degree at a time),
\begin{equation}\label{eq:coloring-kernel}
S_O(x)=\mathcal E\prod_{a\ge1}\det(I+x_aT).
\end{equation}
Consequently, for every $\mu\vdash n$,
\begin{equation}\label{eq:fos-squarefree}
[h_\mu]H_O=\mathcal E\,m_\mu[T].
\end{equation}
\end{proposition}
\begin{proof}
In the factor for color $a$, the term $e_k[T]x_a^k$ selects a
chain of $k$ matching edges to receive that color.  By
\Cref{lem:chain-minors}, its weight records the exact opener and closer
sets, and there is no additional ordering factor.  Squarefree extraction
requires that the selected chains use every opener and every closer
exactly once.  Their union is therefore a perfect matching with opening
set $O$.  Since independent sets of $\Inc(P_M)$ are exactly chains of
$P_M$, these choices are precisely the proper colorings counted by
$S_O$.  This proves \eqref{eq:coloring-kernel}, first with finitely many
colors and then by stabilization.

Apply the matrix evaluation in the second alphabet of the dual Cauchy
identity \cite[Chapter~7]{StanleyEC2}:
\begin{equation}\label{eq:dual-Cauchy-matrix}
\prod_{a\ge1}\det(I+x_aT)
=\sum_\lambda e_\lambda(x)m_\lambda[T].
\end{equation}
Only partitions of $n$ survive $\mathcal E$, because $m_\lambda[T]$ has
degree $|\lambda|$ in each of $u$ and $v$.  Hence
$S_O=\sum_{\mu\vdash n}e_\mu\,\mathcal E m_\mu[T]$.
Applying $\omega$ yields \eqref{eq:fos-squarefree}.
\end{proof}

\section{Coefficient transfer and fixed-opening-set \texorpdfstring{$h$}{h}-positivity}\label{sec:hpositive}

We now convert the squarefree matrix extraction into an elementary coefficient of an auxiliary graph.

\subsection{The auxiliary graph and the coefficient-transfer theorem}\label{subsec:auxiliary-graph}

Retain the threshold data $b_t,B_j,D_O$ and $\eta_j$ from \cref{subsec:thresholds}. Define
\begin{equation}\label{eq:fos-gamma}
 \gamma_{2j-1}=n-\eta_j\quad(1\le j\le n),
 \qquad \gamma_{2j}=\eta_j-1\quad(1\le j<n).
\end{equation}
All entries are nonnegative, and $\eta_n=1$ gives
\begin{equation}\label{eq:fos-gamma-degree}
 |\gamma|=\sum_{s=1}^{2n-1}\gamma_s
 =\sum_{j=1}^{n-1}\bigl((n-\eta_j)+(\eta_j-1)\bigr)+(n-1)
 =n(n-1).
\end{equation}
Write $\gamma!=\prod_{s=1}^{2n-1}\gamma_s!$.

Let $Q$ be the path on vertex set $[2n-1]$, with edges between consecutive integers. Define its
clique blow-up $G_O=Q^\gamma$ by
\begin{align}
 V(G_O)&=\{(s,a):1\le s\le2n-1,\ 1\le a\le\gamma_s\},
 \label{eq:fos-blowup-vertices}\\
 \{(s,a),(t,b)\}\in E(G_O)
 &\quad\Longleftrightarrow\quad |s-t|\le1
 \quad\text{for distinct vertices}.
 \label{eq:fos-blowup-adjacency}
\end{align}
If $\gamma_s=0$, position $s$ is empty: its two neighboring cliques are \emph{not} joined across that gap.

For $\mu\vdash n$, pad $\mu$ with zeros to length $n$ and put
\begin{equation}\label{eq:fos-rectangle-complement}
 \nu_\mu=(n-\mu_n,n-\mu_{n-1},\ldots,n-\mu_1),
\end{equation}
omitting zero parts. This is a partition of $n(n-1)$.

\begin{theorem}[Coefficient transfer]\label{thm:coefficient-transfer}
For every valid opening set $O\subseteq[2n]$ and every $\mu\vdash n$,
\begin{equation}\label{eq:fos-coefficient-transfer}
 [h_\mu]H_O=\frac{D_O}{\gamma!}[e_{\nu_\mu}]X_{G_O}.
\end{equation}
Moreover, $[h_\mu]H_O$ is an integer divisible by $D_O$.
\end{theorem}

\begin{example}\label{ex:coefficient-transfer}
For $n=4$ and $O=\{1,2,4,6\}$, \cref{ex:running-thresholds} gives
\[
 r=(1,1,2,3),\qquad \eta(O)=(2,2,2,1),\qquad
 \gamma=(2,1,2,1,2,1,3).
\]
Thus $D_O=2$ and $\gamma!=48$. For $\mu=(2,1,1)$, padded to $(2,1,1,0)$, we have
$\nu_\mu=(4,3,3,2)$, so \cref{thm:coefficient-transfer} specializes to
\[
 [h_{(2,1,1)}]H_O=\frac1{24}[e_{(4,3,3,2)}]X_{G_O}.
\]
Here $G_O$ is the clique blow-up of a seven-vertex path with clique sizes
$2,1,2,1,2,1,3$, and it has twelve vertices.
\end{example}

\subsection{Repeated thresholds and a universal matrix}\label{subsec:repeated-thresholds}

The variables below collect openers having the same threshold. Introduce new independent variables $z_1,\ldots,z_n$ and the universal matrices
\begin{equation}\label{eq:fos-universal-matrix}
 L_{ij}=\ind[j\le i],\qquad
 Z=\diag(z_1,\ldots,z_n),\qquad K=LZL^{\mathsf T}V.
\end{equation}

\begin{lemma}[Repeated-threshold extraction]\label{lem:fos-threshold-extraction}
For every $\mu\vdash n$,
\begin{equation}\label{eq:fos-threshold-extraction}
 [h_\mu]H_O
 =D_O[z_1^{b_1}\cdots z_n^{b_n}v_1\cdots v_n]m_\mu[K].
\end{equation}
\end{lemma}
\begin{proof}
Under the polynomial substitution
\begin{equation}\label{eq:fos-threshold-substitution}
 z_t\longmapsto\sum_{i:r_i=t}u_i,
\end{equation}
the matrices $LZL^{\mathsf T}$ and $R^{\mathsf T}UR$ agree: their $(a,b)$ entries both become
\[
 \sum_{i:r_i\le\min(a,b)}u_i.
\]
Cyclic invariance of the trace gives
\[
 \tr\bigl((URVR^{\mathsf T})^k\bigr)
 =\tr\bigl((R^{\mathsf T}URV)^k\bigr)\qquad(k\ge1).
\]
As the power sums generate $\Lambda_{\mathbb Q}$, it follows that
\begin{equation}\label{eq:fos-cyclic-evaluation}
 f[URVR^{\mathsf T}]=f[R^{\mathsf T}URV]
 \qquad(f\in\Lambda_{\mathbb Q}).
\end{equation}
The variable groups in \eqref{eq:fos-threshold-substitution} are disjoint, and
\[
 \left[\prod_{i:r_i=t}u_i\right]
 \left(\sum_{i:r_i=t}u_i\right)^{b_t}=b_t!.
\]
Consequently, for any polynomial $F$ in the $z$-variables, with coefficients possibly depending on $v$,
\begin{equation}\label{eq:fos-group-factorial}
 [u_1\cdots u_n]F\left(\sum_{r_i=1}u_i,\ldots,\sum_{r_i=n}u_i\right)
 =D_O[z_1^{b_1}\cdots z_n^{b_n}]F(z_1,\ldots,z_n).
\end{equation}
For $b_t=0$, the substitution sends $z_t$ to $0$ and the right side takes its constant coefficient,
so the formula still holds. In the running example the substitution is
$(z_1,z_2,z_3,z_4)\mapsto(u_1+u_2,u_3,u_4,0)$, and the factor $D_O=2$ comes from
$[u_1u_2](u_1+u_2)^2=2$. Apply \eqref{eq:fos-cyclic-evaluation} and
\eqref{eq:fos-group-factorial} to \cref{prop:fos-squarefree}.
\end{proof}

\begin{remark}\label{rem:universal-inverse}
All inverses below are taken \emph{before} the threshold substitution, in the universal
Laurent-polynomial ring
\[
 \mathbb Q[z_1^{\pm1},\ldots,z_n^{\pm1},v_1^{\pm1},\ldots,v_n^{\pm1}].
\]
Indeed,
\begin{equation}\label{eq:fos-det-K}
 \det K=\prod_{i=1}^n z_iv_i.
\end{equation}
A specialization with $b_t=0$ can make $R^{\mathsf T}UR$ singular; no inverse of that specialized
matrix is used.
\end{remark}

\subsection{The inverse matrix and a path}\label{subsec:inverse-path}

For a finite simple graph $F$, let $\mathcal I(F)$ be the family of its
independent vertex sets, including the empty set, and let $\mathcal I_k(F)$
consist of those of cardinality $k$. Define the homomorphism
\[
 \Phi_F:\Lambda_{\mathbb Z}\longrightarrow\mathbb Z[t_v:v\in V(F)]
\]
by
\begin{equation}\label{eq:fos-graph-homomorphism}
 \Phi_F(e_k)=\sum_{A\in\mathcal I_k(F)}\prod_{a\in A}t_a.
\end{equation}
This is Stanley's graph analogue of symmetric-function evaluation \cite[Section~2]{StanleyGraphAnalogues}. It is well-defined because the elementary functions freely generate $\Lambda_{\mathbb Z}$. The same notation denotes its extension over $\mathbb Q$.

\begin{lemma}[Path evaluation of the inverse]\label{lem:fos-inverse-path}
Define Laurent monomials
\begin{equation}\label{eq:fos-theta}
 \theta_{2i-1}=(z_iv_i)^{-1}\quad(1\le i\le n),\qquad
 \theta_{2i}=(z_{i+1}v_i)^{-1}\quad(1\le i<n).
\end{equation}
Let $Q$ be the path with vertex set $[2n-1]$ and edge set
$\{\{s,s+1\}:1\le s\le2n-2\}$. Then, for every $f\in\Lambda_{\mathbb Q}$,
\begin{equation}\label{eq:fos-inverse-path}
 f[K^{-1}]=\Phi_Q(f)(\theta_1,\ldots,\theta_{2n-1}).
\end{equation}
\end{lemma}
\begin{proof}
The inverse $J=L^{-1}$ is lower bidiagonal:
\[
 J=\begin{pmatrix}
 1&0&0&\cdots\\
 -1&1&0&\cdots\\
 0&-1&1&\cdots\\
 \vdots&&\ddots&\ddots
 \end{pmatrix},\qquad
 K^{-1}=V^{-1}J^{\mathsf T}Z^{-1}J.
\]
The principal-minor expansion and Cauchy--Binet identity imply, for $0\le k\le n$,
\begin{equation}\label{eq:fos-cauchy-binet}
 e_k[K^{-1}]
 =\sum_{\substack{I,A\subseteq[n]\\|I|=|A|=k}}
     \det J[A,I]^2\prod_{i\in I}v_i^{-1}\prod_{a\in A}z_a^{-1}.
\end{equation}
Rows and columns of each submatrix are taken in increasing order.

The bipartite support graph of $J$ has a vertex for each row and column,
with an edge for each nonzero entry. It is the path
\[
 \mathrm{row}_1-\mathrm{col}_1-\mathrm{row}_2-
 \mathrm{col}_2-\cdots-\mathrm{row}_n-\mathrm{col}_n.
\]
This support path has $2n$ vertices and $2n-1$ edges.
A nonzero term in the determinant expansion of $J[A,I]$ corresponds to a
perfect matching of its bipartite support graph: it selects exactly one
nonzero entry in each row indexed by $A$ and each column indexed by $I$.
There is at most one such matching, since the symmetric difference of two
distinct perfect matchings would contain a cycle, whereas the support is
a subgraph of a path. The nonzero entries of $J[A,I]$ are $\pm1$, so
$\det J[A,I]^2$ is $1$ exactly when such a perfect matching exists, and is
$0$ otherwise.

It follows that \eqref{eq:fos-cauchy-binet} counts every size-$k$ matching
of the support path once, weighted by the inverses of the variables at its
endpoints. Label its edges $1,\ldots,2n-1$ in the displayed order: the
diagonal edge $\{\mathrm{row}_i,\mathrm{col}_i\}$ has label $2i-1$, and the
subdiagonal edge $\{\mathrm{row}_{i+1},\mathrm{col}_i\}$ has label $2i$.
Their weights are
\[
 (z_iv_i)^{-1}=\theta_{2i-1}\quad(1\le i\le n),\qquad
 (z_{i+1}v_i)^{-1}=\theta_{2i}\quad(1\le i<n).
\]
The line graph has one vertex for each of these $2n-1$ edges, with two
vertices adjacent exactly when their labels are consecutive. Thus it is
$Q$, which has $2n-1$ vertices and $2n-2$ edges. Matchings of the support
path correspond exactly to independent sets in $Q$. Hence
$e_k[K^{-1}]=\Phi_Q(e_k)(\theta)$. For $k>n$ both sides vanish, and for
$k=0$ both equal $1$. Since the $e_k$ generate $\Lambda_{\mathbb Q}$, the identity holds for every $f$.
\end{proof}

\begin{lemma}[Clique-blow-up coefficients]\label{lem:fos-clique-blowup}
Let $F$ be a finite simple graph and $\alpha=(\alpha_v)_{v\in V(F)}$ a nonnegative integer vector. Let $F^\alpha$ replace each vertex $v$ by a clique of size $\alpha_v$, with all edges between two such cliques when their original vertices are adjacent. For every $\lambda\vdash\sum_v\alpha_v$,
\begin{equation}\label{eq:fos-clique-blowup}
 [t^\alpha]\Phi_F(m_\lambda)
   =\frac{1}{\prod_v\alpha_v!}[e_\lambda]X_{F^\alpha},
 \qquad t^\alpha:=\prod_v t_v^{\alpha_v}.
\end{equation}
\end{lemma}
\begin{proof}
This is the coefficient form of the graph-analogue construction in \cite[Proposition~2.1 and Corollary~2.2(b)]{StanleyGraphAnalogues}; we include the argument to fix the factorial. The dual Cauchy identity, in its elementary--monomial form, is
\[
 \prod_{a,b}(1+x_ay_b)=\sum_\lambda e_\lambda(x)m_\lambda(y).
\]
Applying $\Phi_F$ in the $y$-variables gives
\begin{equation}\label{eq:fos-color-set-generating-function}
 \sum_\lambda e_\lambda(x)\Phi_F(m_\lambda)(t)
 =\prod_{a\ge1}\left(\sum_{S\in\mathcal I(F)}
                   x_a^{|S|}\prod_{v\in S}t_v\right).
\end{equation}
The equality is understood degree by degree, or first with finitely many $x$-variables and then stably.

In the factor indexed by color $a$, choose the independent set of original vertices receiving that color. Taking $[t^\alpha]$ therefore assigns to every vertex $v$ a set of exactly $\alpha_v$ distinct colors, with disjoint color sets at adjacent vertices. Each assignment gives exactly $\prod_v\alpha_v!$ proper colorings of the labeled graph $F^\alpha$, by choosing, for each $v$, a bijection from its assigned color set to the vertices in the clique replacing $v$. Conversely, forgetting labels within each clique recovers the assignment. This correspondence preserves the monomial weight in the $x$-variables. Thus
\[
 X_{F^\alpha}(x)
 =\left(\prod_v\alpha_v!\right)[t^\alpha]
          \sum_\lambda e_\lambda(x)\Phi_F(m_\lambda)(t).
\]
Comparing elementary coefficients proves \eqref{eq:fos-clique-blowup}, including the case of zero clique sizes.
\end{proof}

\subsection{Rectangle complementation and the final extraction}\label{subsec:rectangle-complement}

\begin{proof}[Proof of \cref{thm:coefficient-transfer}]
Fix $\mu\vdash n$. In $n$ invertible variables, the rectangular-complement identity is
\begin{equation}\label{eq:fos-complement-identity}
 m_\mu(y_1,\ldots,y_n)
 =(y_1\cdots y_n)^n
        m_{\nu_\mu}(y_1^{-1},\ldots,y_n^{-1}).
\end{equation}
Indeed, the map $(a_1,\ldots,a_n)\mapsto(n-a_1,\ldots,n-a_n)$ bijects the distinct rearrangements of the zero-padded $\mu$ with those of $\nu_\mu$. There is no additional factorial, even when parts repeat.

The identity may be evaluated on the universal matrix $K$. Indeed,
\[
 \mathbb Q[y_1^{\pm1},\ldots,y_n^{\pm1}]^{\mathfrak S_n}
 =\mathbb Q[e_1,\ldots,e_n,e_n^{-1}],
\]
where the $e_i$ on the right are evaluated in $y_1,\ldots,y_n$.
To see this, multiply any symmetric Laurent polynomial by a sufficiently
large power of $e_n=y_1\cdots y_n$ to clear its negative exponents. The
result is a symmetric polynomial, hence a polynomial in $e_1,\ldots,e_n$.
Thus only $e_n^{-1}$ needs to be adjoined. Since $e_n[K]=\det K$ is a unit
in the universal Laurent ring of \cref{rem:universal-inverse}, matrix
evaluation extends to this ring of symmetric Laurent polynomials.
Moreover,
\[
 \det(I+tK^{-1})=\det(K)^{-1}\det(K+tI)
\]
implies
\[
 e_k[K^{-1}]=\frac{e_{n-k}[K]}{e_n[K]}\qquad(0\le k\le n).
\]
Thus matrix inversion implements finite-variable inversion in \eqref{eq:fos-complement-identity}. Put $\Omega=\prod_{i=1}^n z_iv_i=\det K$. By \cref{lem:fos-inverse-path},
\begin{equation}\label{eq:fos-complement-path}
 m_\mu[K]=\Omega^n\Phi_Q(m_{\nu_\mu})(\theta).
\end{equation}
This step takes place in the universal Laurent ring of \cref{rem:universal-inverse}.

Since $\Phi_Q$ preserves degree, $\Phi_Q(m_{\nu_\mu})(t)$ is homogeneous
of total degree $|\nu_\mu|=n(n-1)$. Write its finite monomial expansion as
\[
 \Phi_Q(m_{\nu_\mu})(t)
 =\sum_{\substack{\delta\in\mathbb Z_{\ge0}^{2n-1}\\|\delta|=n(n-1)}}
       a_\delta t^\delta,
 \qquad t^\delta:=\prod_{s=1}^{2n-1}t_s^{\delta_s}.
\]
Here $\delta=(\delta_1,\ldots,\delta_{2n-1})$ is an exponent vector and
$|\delta|:=\sum_{s=1}^{2n-1}\delta_s$. Set the two auxiliary boundary
values $\delta_0=\delta_{2n}=0$. By
\eqref{eq:fos-theta}, the exponents of $z_i$ and $v_i$ in
$\Omega^n\theta^\delta$ are, respectively,
\begin{equation}\label{eq:fos-exponents}
 n-\delta_{2i-2}-\delta_{2i-1},\qquad
 n-\delta_{2i-1}-\delta_{2i}.
\end{equation}
Consequently a term contributes to $z_1^{b_1}\cdots z_n^{b_n}v_1\cdots v_n$ if and only if
\begin{equation}\label{eq:fos-exponent-system}
 \begin{aligned}
 \delta_{2i-2}+\delta_{2i-1}&=n-b_i,\\
 \delta_{2i-1}+\delta_{2i}&=n-1
 \end{aligned}
 \qquad(1\le i\le n).
\end{equation}
Subtracting the first equation from the second gives
$\delta_{2i}-\delta_{2i-2}=b_i-1$. Starting from $\delta_0=0$, we obtain uniquely
\begin{align}
 \delta_{2i}&=\sum_{t=1}^i b_t-i=B_i-i=\eta_i-1,
                  \label{eq:fos-even-exponents}\\
 \delta_{2i-1}&=n-1-\delta_{2i}=n-\eta_i.
                  \label{eq:fos-odd-exponents}
\end{align}
The final boundary condition holds because $\delta_{2n}=B_n-n=0$. The entries in \eqref{eq:fos-even-exponents}--\eqref{eq:fos-odd-exponents} are nonnegative by \eqref{eq:eta-basic}, and they are exactly $\gamma$. Thus no other monomial of $\Phi_Q(m_{\nu_\mu})$ contributes.

Combining \cref{lem:fos-threshold-extraction} with \eqref{eq:fos-complement-path} yields the integral coefficient formula
\begin{equation}\label{eq:fos-integral-coefficient}
 [h_\mu]H_O=D_O[t^\gamma]\Phi_Q(m_{\nu_\mu}).
\end{equation}
Apply \cref{lem:fos-clique-blowup} with $F=Q$, $\alpha=\gamma$, and $\lambda=\nu_\mu$. We obtain
\[
 [h_\mu]H_O=\frac{D_O}{\gamma!}[e_{\nu_\mu}]X_{Q^\gamma},
\]
which is \eqref{eq:fos-coefficient-transfer}.

Since $m_{\nu_\mu}\in\Lambda_{\mathbb Z}$ and $\Phi_Q$ is defined over
$\mathbb Z$, formula \eqref{eq:fos-integral-coefficient} also shows that
$[h_\mu]H_O$ is an integer divisible by $D_O$.
The calculation includes $n=1$: then $\gamma=(0)$,
$\nu_{(1)}=\varnothing$, and $X_\varnothing=e_\varnothing=1$.
\end{proof}

\subsection{The unit interval graph and the positivity input}\label{subsec:unit-interval}

A unit interval graph is the intersection graph of intervals of equal length: vertices are
adjacent exactly when their intervals intersect.

\begin{theorem}[Fixed-opening-set positivity]\label{thm:fos-h-positive}
For every valid opening set $O\subseteq[2n]$,
\begin{equation}\label{eq:fos-positivity}
 H_O=\sum_{\mu\vdash n}c^O_\mu h_\mu,
 \qquad S_O=\sum_{\mu\vdash n}c^O_\mu e_\mu,
 \qquad c^O_\mu\in\mathbb Z_{\ge0}.
\end{equation}
\end{theorem}
\begin{proof}
For $n=1$, directly $H_{\{1\}}=h_1$; assume $n\ge2$.
Associate to every vertex in the clique at position $s$ the unit interval
\begin{equation}\label{eq:fos-unit-intervals}
 I_s=\left[\frac{2s}{3},\frac{2s}{3}+1\right].
\end{equation}
Intervals from the same or consecutive positions intersect, whereas intervals from positions at
distance at least two are disjoint. They realize exactly \eqref{eq:fos-blowup-adjacency}. When a
position is empty, no interval is assigned there and the original spacing is unchanged.

To match Hikita's natural-labeling convention directly, label the vertices first by increasing
clique position and then in any fixed order within each clique, and put
\[
 g_s=\sum_{t=1}^s\gamma_t\quad(s\ge1),\qquad g_s=0\quad(s\le0).
\]
If label $i$ is in the clique at position $s$, define
\begin{equation}\label{eq:fos-Hessenberg}
 \mathsf e(i)=g_{s-2}.
\end{equation}
Then $0\le\mathsf e(i)<i$, the sequence $\mathsf e(i)$ is weakly increasing, and the earlier
neighbors of $i$ are exactly
$\mathsf e(i)+1,\ldots,i-1$. Hence this is a conjugate Hessenberg function and $G_O$ is its
naturally labeled unit interval graph in the sense of \cite[Definition~2.3]{Hikita2025SS}.

Hikita's Corollary~1.7 in arXiv version~2 \cite{Hikita2025SS}, specialized to $q=1$, gives
\[
 [e_{\nu_\mu}]X_{G_O}\ge0.
\]
Only the ordinary chromatic symmetric function is used: at $q=1$ the ascent weight is $1$ for every
proper coloring, so Hikita's chromatic quasisymmetric function specializes to $X_{G_O}$. Since
$D_O/\gamma!>0$, the transfer formula proves $[h_\mu]H_O\ge0$.

Integrality follows from \cref{thm:coefficient-transfer}.
Applying $\omega(h_\mu)=e_\mu$ gives the stated expansion of $S_O$.
\end{proof}

\begin{remark}\label{rem:scope-positivity}
Hikita's result is applied to the auxiliary graph $G_O$, not to the individual matching graphs
$\Inc(P_M)$. The argument uses only $q=1$ and gives no coefficientwise assertion in $q$.
\end{remark}

\begin{corollary}[Integrality after threshold regrouping]\label{cor:Fr-positive}
For every valid threshold sequence $r$, the symmetric function $F_r$ of \cref{prop:regrouping} is
$e$-positive with integer coefficients. Moreover, for every $\mu\vdash n$,
\begin{equation}\label{eq:Fr-coefficient}
 [e_\mu]F_r=[t^\gamma]\Phi_Q(m_{\nu_\mu})\in\mathbb Z_{\ge0},
 \qquad D_O\mid [h_\mu]H_O.
\end{equation}
\end{corollary}
\begin{proof}
Combine $S_O=D_OF_r$, \eqref{eq:h-e-duality}, and \eqref{eq:fos-integral-coefficient}; positivity
follows from \cref{thm:fos-h-positive}.
\end{proof}

\section{Power-sum refinements}\label{sec:autocorrelation}

\subsection{A cycle-type autocorrelation identity}
Recall from \cref{prop:threshold} that $P_r$ indexes the matchings in $\Mcal(O)$. For
$\lambda\vdash n$, put
\begin{equation}\label{eq:pair-count}
A_r(\lambda):=
\#\{(\pi,\rho)\in P_r^2:\cyc(\rho^{-1}\pi)=\lambda\}.
\end{equation}

\begin{lemma}\label{lem:trace-count}
For every $\lambda\vdash n$,
\begin{equation}\label{eq:trace-count}
\mathcal E\,p_\lambda[T]=z_\lambda A_r(\lambda).
\end{equation}
\end{lemma}
\begin{proof}
Fix a permutation $\sigma\in \mathfrak S_n$ of cycle type $\lambda$.
Expanding the traces along its cycles gives
\begin{equation}\label{eq:fos-trace-expansion}
p_\lambda[T]
=\sum_{\alpha,\beta:[n]\to[n]}
\prod_{x=1}^n
u_{\alpha(x)}v_{\beta(x)}
R_{\alpha(x),\beta(x)}R_{\alpha(\sigma(x)),\beta(x)}.
\end{equation}
A term contributes to $\mathcal E$ only when both $\alpha$ and $\beta$
are bijections.  For such a nonzero term, define
\[
\pi:=\beta\alpha^{-1},\qquad
\rho:=\beta\sigma^{-1}\alpha^{-1}.
\]
The first and second $R$-factors in \eqref{eq:fos-trace-expansion} imply
$\pi\in P_r$ and $\rho\in P_r$, respectively. Moreover,
\begin{equation}\label{eq:fos-relative-permutation}
\rho^{-1}\pi=\alpha\sigma\alpha^{-1}.
\end{equation}
Conversely, fix $(\pi,\rho)\in P_r^2$ of the required relative cycle type.
There are exactly $z_\lambda$ bijections $\alpha$ satisfying
\eqref{eq:fos-relative-permutation}: the set of such conjugators $\alpha$ is a coset
of the centralizer of $\sigma$, whose order is
$z_\lambda=\prod_j j^{m_j(\lambda)}m_j(\lambda)!$ (rotate cycles and
permute cycles of equal length). For each such $\alpha$, $\beta=\pi\alpha$ is uniquely
determined and gives a nonzero squarefree term.  The lemma follows.
\end{proof}

For the indicator function $f_r=\mathbf 1_{P_r}$ on $\mathfrak S_n$,
the function
\[
C_r(\tau):=\sum_{\rho\in\mathfrak S_n}f_r(\rho)f_r(\rho\tau)
\]
is the (right) group autocorrelation of $f_r$. In particular,
$A_r(\lambda)=\sum_{\cyc(\tau)=\lambda}C_r(\tau)$.

\begin{theorem}[Cycle-type autocorrelation identity]\label{thm:autocorrelation}
For every valid opening set $O$,
\begin{equation}\label{eq:fos-autocorrelation}
H_O=\sum_{\pi,\rho\in P_r}p_{\cyc(\rho^{-1}\pi)}
=\sum_{\lambda\vdash n}A_r(\lambda)p_\lambda.
\end{equation}
\end{theorem}
\begin{proof}
By \Cref{prop:fos-squarefree} and Hall duality,
$\langle H_O,m_\mu\rangle=\mathcal E m_\mu[T]$ for every $\mu\vdash n$.
Since the monomial symmetric functions form a basis, it follows that
$\langle H_O,f\rangle=\mathcal E f[T]$ for every homogeneous symmetric function $f$ of degree $n$.
Taking $f=p_\lambda$ and using \Cref{lem:trace-count} gives
$\langle H_O,p_\lambda\rangle=z_\lambda A_r(\lambda)$.
Since $\langle p_\lambda,p_\mu\rangle=\delta_{\lambda\mu}z_\lambda$,
the coefficient of $p_\lambda$ in $H_O$ is $A_r(\lambda)$.
\end{proof}

\subsection{Explicit power-sum coefficients}\label{subsec:power-sum-coefficients}

The autocorrelation identity expresses each power-sum coefficient of $H_O$ through pairs of restricted
permutations. The next formula depends only on the threshold sequence $r$.

For $\lambda\vdash n$, let
\[
\mathfrak C_\lambda
=
\{\tau\in\mathfrak S_n:\cyc(\tau)=\lambda\}
\]
be the conjugacy class of cycle type $\lambda$. In view of \cref{thm:autocorrelation}, define
\begin{equation}\label{eq:a-lambda-O}
a_{\lambda,O}
:=
[p_\lambda]H_O
=
A_r(\lambda).
\end{equation}

We use the following count for permutations with position-dependent lower bounds.

\begin{lemma}[Nested lower-bound count]\label{lem:nested-lower-bound}
Let $d_1,\ldots,d_n$ be integers in $[n]$, and let
\[
d_{(1)}\le d_{(2)}\le\cdots\le d_{(n)}
\]
be their weakly increasing rearrangement. Then
\begin{equation}\label{eq:nested-lower-bound}
\#\{\sigma\in\mathfrak S_n:\sigma(j)\ge d_j\text{ for every }j\}
=
\begin{cases}
\displaystyle\prod_{i=1}^n\bigl(i-d_{(i)}+1\bigr),
& d_{(i)}\le i\text{ for all }i,\\[6pt]
0,&\text{otherwise}.
\end{cases}
\end{equation}
\end{lemma}

\begin{proof}
Order the positions as $j_1,\ldots,j_n$ so that $d_{j_i}=d_{(i)}$, and assign their values in the
order $j_n,\ldots,j_1$. When position $j_i$ is reached, all $n-i$ previously assigned values lie in
$\{d_{(i)},\ldots,n\}$. There are therefore
\[
 n-d_{(i)}+1-(n-i)=i-d_{(i)}+1
\]
available values. Multiplication gives the formula when $d_{(i)}\le i$ for every $i$. If
$d_{(i)}>i$, the $n-i+1$ positions $j_i,\ldots,j_n$ require distinct values from a set of size
$n-d_{(i)}+1<n-i+1$, so no such permutation exists.
\end{proof}

For $\tau\in\mathfrak S_n$, define
\begin{equation}\label{eq:s-tau}
s_j(\tau)
:=
\max\{r_j,r_{\tau^{-1}(j)}\},
\qquad
1\le j\le n,
\end{equation}
and let
\[
s_{(1)}(\tau)\le\cdots\le s_{(n)}(\tau)
\]
be the weakly increasing rearrangement of $s_1(\tau),\ldots,s_n(\tau)$. Set
\begin{equation}\label{eq:tau-weight}
\operatorname{wt}_r(\tau)
:=
\begin{cases}
\displaystyle\prod_{i=1}^n
\bigl(i-s_{(i)}(\tau)+1\bigr),
& s_{(i)}(\tau)\le i\text{ for all }i,\\[6pt]
0,&\text{otherwise}.
\end{cases}
\end{equation}

\begin{proposition}[Explicit power-sum coefficients]
\label{thm:explicit-power-sum-coefficient}
For every valid opening set $O\subseteq[2n]$ and every partition $\lambda\vdash n$,
\begin{equation}\label{eq:explicit-power-sum-coefficient}
a_{\lambda,O}
=
\sum_{\tau\in\mathfrak C_\lambda}
\operatorname{wt}_r(\tau).
\end{equation}
Equivalently,
\[
H_O
=
\sum_{\lambda\vdash n}
\left(
\sum_{\tau\in\mathfrak C_\lambda}
\operatorname{wt}_r(\tau)
\right)p_\lambda.
\]
\end{proposition}

\begin{proof}
Fix $\tau\in\mathfrak S_n$. The relation $\rho^{-1}\pi=\tau$ determines $\pi=\rho\tau$. Both
permutations belong to $P_r$ precisely when
\[
 \rho(j)\ge r_j,\qquad \rho(j)\ge r_{\tau^{-1}(j)}
 \qquad(j\in[n]),
\]
or equivalently $\rho(j)\ge s_j(\tau)$. By \cref{lem:nested-lower-bound}, there are
$\operatorname{wt}_r(\tau)$ choices for $\rho$. Summing over $\tau\in\mathfrak C_\lambda$ in
\cref{thm:autocorrelation} proves the formula.
\end{proof}

\subsection{Alternating permutations and record partitions}
\label{subsec:fixed-valley-power-sum}

The coefficients $a_{\lambda,O}$ refine the global record-partition interpretation in \cite[Theorem~3.2]{ASS}.

Let
\[
w=a_1a_2\cdots a_{2n}\in\mathfrak S_{2n}
\]
be alternating, that is,
\begin{equation}\label{eq:alternating-convention}
a_1>a_2<a_3>a_4<\cdots>a_{2n-1}>a_{2n}.
\end{equation}
Its valley-value set is
\[
\Val(w):=\{a_2,a_4,\ldots,a_{2n}\},
\]
and its odd-position word is
\[
\widehat w:=a_1a_3\cdots a_{2n-1}.
\]

\begin{definition}[Record partition]\label{def:record-partition}
A position $j$ in a list $b_1,\ldots,b_m$ is a left-to-right record if $b_j>b_i$ for every $i<j$.
The first position is always a record. Let the left-to-right record positions of $\widehat w$ be
\[
1=q_1<q_2<\cdots<q_s.
\]
The \emph{record composition} of $\widehat w$ is
\[
(q_2-q_1,\ q_3-q_2,\ \ldots,\ q_s-q_{s-1},\ n+1-q_s),
\]
and the \emph{record partition}, denoted $\rp(\widehat w)$, is the decreasing rearrangement of these parts.
\end{definition}

This agrees with the record-partition convention of \cite[Section~3]{ASS}. For a valid opening set
$O\subseteq[2n]$, define
\begin{equation}\label{eq:Alt-O}
\Alt_{2n}(O)
:=
\{w\in\mathfrak S_{2n}:w\text{ satisfies \eqref{eq:alternating-convention} and }\Val(w)=O\}.
\end{equation}

\begin{example}\label{ex:running-records}
For $O=\{1,2,4,6\}$, the alternating permutation $w=31745286$ has
$\Val(w)=\{1,4,2,6\}=O$ and $\widehat w=3758$. The record positions are $1,2,4$,
so the record composition is $(1,2,1)$ and $\rp(\widehat w)=(2,1,1)$.
\end{example}

\begin{proposition}[Bicolored matching correspondence]
\label{prop:bicolored-alternating}
There is a bijection
\[
\Psi_O:
\Mcal(O)\times\Mcal(O)
\longrightarrow
\Alt_{2n}(O)
\]
such that, if $M_\pi$ and $M_\rho$ correspond to $\pi,\rho\in P_r$ under \cref{prop:threshold}
and $w:=\Psi_O(M_\pi,M_\rho)$, then
\begin{equation}\label{eq:record-cycle-correspondence}
\rp(\widehat w)=\cyc(\rho^{-1}\pi).
\end{equation}
\end{proposition}

\begin{proof}
Take an ordered pair
\[
(M_\pi,M_\rho)\in\Mcal(O)\times\Mcal(O).
\]
Color the edges of $M_\pi$ red and the edges of $M_\rho$ blue. If the two matchings share an edge,
we retain two copies of that edge, one of each color. Thus the red-blue union is a $2$-regular
edge-colored multigraph on the vertex set $[2n]$. Every connected component is an even cycle whose
edge colors alternate; a common edge forms the two-edge component consisting of its red and blue copies.

Every edge of either matching joins an opener in $O$ to a closer in $C=[2n]\setminus O$. In
particular, every opener is smaller than both of its neighbors in the red-blue union. It follows
that the largest vertex of each connected component is a closer.

Fix a connected component $C_0$, and let $p_1$ be its largest vertex. Starting at $p_1$, traverse
first the red edge incident with $p_1$ and then alternate blue and red edges until returning to
$p_1$. Because the colors alternate around the component, this traversal gives a word
\begin{equation}\label{eq:component-word}
B_{C_0}
=
p_1v_1p_2v_2\cdots p_kv_k,
\end{equation}
where the $v_i$ are openers and the $p_i$ are closers. Each matching edge has its opener as its
smaller endpoint, so
\[
p_i>v_i
\qquad(1\le i\le k),
\]
and, because the edge from $v_i$ to $p_{i+1}$ is blue for $i<k$,
\[
v_i<p_{i+1}
\qquad(1\le i<k).
\]
The final blue edge joins $v_k$ to $p_1$, so $v_k<p_1$ as well. Therefore \eqref{eq:component-word}
is an alternating word beginning with the largest vertex of its component.

Let the connected components be
\[
C_1,\ldots,C_s,
\]
ordered so that their largest vertices satisfy
\[
p_1^{(1)}<p_1^{(2)}<\cdots<p_1^{(s)}.
\]
Define
\begin{equation}\label{eq:Psi-forward}
\Psi_O(M_\pi,M_\rho)
:=
B_{C_1}B_{C_2}\cdots B_{C_s}.
\end{equation}
The alternating inequalities hold inside every block. At the boundary between $B_{C_j}$ and
$B_{C_{j+1}}$, the last entry $v$ of $B_{C_j}$ satisfies
\[
v<p_1^{(j)}<p_1^{(j+1)},
\]
so the required valley-to-peak inequality also holds across the boundary. Hence the concatenation is
alternating. Its even-position entries are precisely the openers of all components, so its
valley-value set is exactly $O$. Thus
\[
\Psi_O(M_\pi,M_\rho)\in\Alt_{2n}(O).
\]

We construct the inverse map. Let
\[
w=p_1v_1p_2v_2\cdots p_nv_n\in\Alt_{2n}(O),
\]
where the $p_i$ are the odd-position entries and the $v_i$ are the valley values. Break $w$
immediately before each odd-position entry that is a left-to-right record of
\[
\widehat w=p_1p_2\cdots p_n.
\]
Consider one resulting block
\[
B=p_1v_1p_2v_2\cdots p_kv_k.
\]
By the definition of the record blocks, no $p_i$ with $i>1$ exceeds $p_1$. Each valley $v_i$ is
smaller than the peak immediately preceding it. Hence $p_1$ is the largest entry of the entire
block.

On the vertices of $B$, define a red matching by
\begin{equation}\label{eq:red-inverse}
M_{\mathrm{red}}(B)
=
\{\{v_i,p_i\}:1\le i\le k\},
\end{equation}
and a blue matching by
\begin{equation}\label{eq:blue-inverse}
M_{\mathrm{blue}}(B)
=
\{\{v_i,p_{i+1}\}:1\le i<k\}
\cup
\{\{v_k,p_1\}\}.
\end{equation}
The alternating inequalities give
\[
v_i<p_i
\qquad(1\le i\le k)
\]
and
\[
v_i<p_{i+1}
\qquad(1\le i<k).
\]
The remaining inequality $v_k<p_1$ holds because $p_1$ is the largest entry of the block. Therefore
every $v_i$ is the smaller endpoint of its red edge and of its blue edge.

Applying \eqref{eq:red-inverse}--\eqref{eq:blue-inverse} to all record blocks produces two perfect matchings
\[
M_{\mathrm{red}},M_{\mathrm{blue}}\in\Mcal(O).
\]
Within each block their red-blue union is one alternating cycle, and the first peak is the largest
vertex of that component. Moreover, the first peaks of the blocks are precisely the left-to-right
records of $\widehat w$, so they occur in increasing order. Therefore, when the forward construction
\eqref{eq:Psi-forward} is applied to the two matchings reconstructed from $w$, it traverses each
block in the same order and concatenates the blocks in the same record order; the resulting word is
exactly $w$. Conversely, if one starts from a pair of matchings, forms the component words, and then
cuts the resulting alternating permutation at the record peaks, each cut occurs precisely at the
first vertex of one red-blue component. The inverse construction therefore recovers every red and
blue edge of the original pair. Hence the two maps are mutually inverse, and $\Psi_O$ is a
bijection.

Set $w:=\Psi_O(M_\pi,M_\rho)$. It remains to verify that \eqref{eq:record-cycle-correspondence} holds.
Label the openers as
$o_1<\cdots<o_n$ and the closers as $c_1<\cdots<c_n$. Under \cref{prop:threshold},
\[
M_\pi:\quad o_i\longleftrightarrow c_{\pi(i)},
\qquad
M_\rho:\quad o_i\longleftrightarrow c_{\rho(i)}.
\]
Starting at the opener $o_i$, following its red edge to $c_{\pi(i)}$, and then following the blue
edge incident with that closer back to an opener gives
\[
o_i
\longmapsto
o_{\rho^{-1}\pi(i)}.
\]
Thus the indices of the openers in each red-blue component form one cycle of $\rho^{-1}\pi$,
whose length is the number of openers in that component.

On the alternating-permutation side, the first peak of every component block is larger than every
other entry of that block, while the first peaks of successive blocks are strictly increasing.
Therefore the left-to-right record positions of $\widehat w$ occur exactly at
the beginnings of the component blocks. The lengths of the record blocks in the odd-position word are consequently the
numbers of openers in the red-blue components. Sorting these lengths gives
\[
\rp(\widehat w)
=
\cyc(\rho^{-1}\pi),
\]
which proves \eqref{eq:record-cycle-correspondence}.
\end{proof}

\begin{theorem}[Fixed-opening-set power-sum refinement]
\label{thm:fixed-valley-power-sum}
For every valid opening set $O\subseteq[2n]$,
\begin{equation}\label{eq:fixed-valley-power-sum}
H_O
=
\sum_{w\in\Alt_{2n}(O)}
p_{\rp(\widehat w)}.
\end{equation}
Equivalently, for every $\lambda\vdash n$,
\begin{equation}\label{eq:fixed-valley-power-sum-coeff}
a_{\lambda,O}
=
\#\{w\in\Alt_{2n}(O):\rp(\widehat w)=\lambda\}.
\end{equation}
\end{theorem}

\begin{proof}
By \cref{thm:autocorrelation},
\[
H_O
=
\sum_{\pi,\rho\in P_r}
p_{\cyc(\rho^{-1}\pi)}.
\]
The threshold bijection \cref{prop:threshold} identifies $P_r^2$ with $\Mcal(O)^2$. Applying the
bijection in \cref{prop:bicolored-alternating} and then using \eqref{eq:record-cycle-correspondence}
transforms the preceding sum into
\[
\sum_{w\in\Alt_{2n}(O)}
p_{\rp(\widehat w)}.
\]
This proves \eqref{eq:fixed-valley-power-sum}. Taking the coefficient of $p_\lambda$ gives
\eqref{eq:fixed-valley-power-sum-coeff}.
\end{proof}

\begin{corollary}[Refinement of the global record-partition count]
\label{cor:global-record-refinement}
For every $\lambda\vdash n$, summing over all valid opening sets $O\subseteq[2n]$,
\begin{equation}\label{eq:global-record-refinement}
\sum_O a_{\lambda,O}
=
\#\{w\in\mathfrak S_{2n}:w\text{ is alternating and }\rp(\widehat w)=\lambda\}.
\end{equation}
By \cite[Theorem~3.2]{ASS}, the common value in \eqref{eq:global-record-refinement} is
$|\phi(\lambda)|$ in the notation of that paper.
\end{corollary}

\begin{proof}
Every alternating permutation $w=a_1\cdots a_{2n}$ has a unique valley-value set
\[
O=\{a_2,a_4,\ldots,a_{2n}\}.
\]
This set is valid: pairing each valley $a_{2i}$ with the preceding peak $a_{2i-1}$ gives a perfect
matching whose opening set is $O$, because $a_{2i}<a_{2i-1}$. Hence the sets $\Alt_{2n}(O)$, as $O$
ranges over valid opening sets, form a disjoint partition of all alternating permutations of $[2n]$.
Summing \eqref{eq:fixed-valley-power-sum-coeff} over $O$ proves \eqref{eq:global-record-refinement}.
The final statement is exactly the global record-partition interpretation in
\cite[Theorem~3.2]{ASS}.
\end{proof}

The autocorrelation identity also follows from Chow's path--cycle reciprocity
\cite[Chapter~2, Section~2, Theorem~1]{ChowThesis}. Fix $\pi\in P_r$ and let
$D_\pi$ be the strict comparability digraph of $P_{M_\pi}$, with vertices
indexed by opener ranks. Its arcs are $i\to j$ exactly when $\pi(i)<r_j$.
Its directed complement $D_\pi^c$, including loops, therefore has $i\to j$
exactly when $\pi(i)\ge r_j$. Taking the cycle-cover specialization on the
$D_\pi^c$ side of Chow's reciprocity gives
\[
 \omega X_{\Inc(P_{M_\pi})}
 =
 \sum_{\substack{\sigma\in\mathfrak S_n\\
                  \pi(i)\ge r_{\sigma(i)}\ \text{for all }i}}
 p_{\cyc(\sigma)},
\]
since $D_\pi$ is acyclic and its paths are chains. Each cycle cover is counted
once as a permutation $\sigma$. The condition in the sum is equivalent to
$\rho:=\pi\sigma^{-1}\in P_r$, and then $\rho^{-1}\pi=\sigma$.
Summing over $\pi\in P_r$ recovers \eqref{eq:fos-autocorrelation}.

\section{The fixed-opening-set Schur expansion}\label{sec:schur}

Recall the staggered shape $\rho(\lambda)$ and the coefficient $g_O^\lambda$ from
\eqref{eq:rho-shape} and \eqref{eq:gOlambda}.

\begin{theorem}[Fixed-opening-set Schur expansion]\label{thm:fixedO}
For every valid opening set $O\subseteq[2n]$,
\begin{equation}\label{eq:fixedO-schur}
S_O
=
\sum_{\lambda\vdash n}g_O^\lambda s_{\lambda'}.
\end{equation}
Equivalently,
\begin{equation}\label{eq:fixedO-schur-omega}
H_O
=
\sum_{\lambda\vdash n}g_O^\lambda s_\lambda.
\end{equation}
\end{theorem}

\subsection{\texorpdfstring{$O$}{O}-arrays}\label{sec:oarrays}

\begin{definition}[$O$-array]\label{def:Oarray}
Let $\alpha=(\alpha_1,\ldots,\alpha_\ell)$ be a weak composition of $n$.  An \emph{$O$-array of type
$\alpha$} is an ordered list of rows
\[
R=(R_1,\ldots,R_\ell)
\]
satisfying the following conditions.
\begin{enumerate}
\item Row $R_i$ contains exactly $2\alpha_i$ entries.
\item The entries in every row are strictly increasing from left to right.
\item Every integer in $[2n]$ occurs exactly once in the entire array.
\item If $R_i=(r_{i,1},\ldots,r_{i,2\alpha_i})$, then
\[
r_{i,2j-1}\in O,
\qquad
r_{i,2j}\notin O.
\]
\end{enumerate}
The set of such arrays is denoted by $\Acal_O(\alpha)$.
A zero component of $\alpha$ corresponds to an empty row.
\end{definition}

A row records a chain of matching edges by pairing consecutive entries:
\[
(r_{i,1},r_{i,2}),
(r_{i,3},r_{i,4}),\ldots.
\]
For $O=\{1,2,4,6\}$, the rows $R_1=(1,3,4,7)$ and $R_2=(2,5,6,8)$ from
\cref{ex:running-tableau} form an $O$-array of type $(2,2)$. They record the chains
$13<_{P_M}47$ and $25<_{P_M}68$, and recover the same matching
$M=13\mid25\mid47\mid68$.

\begin{lemma}[Colorings and $O$-arrays]\label{lem:Oarray}
For every weak composition $\alpha=(\alpha_1,\ldots,\alpha_\ell)$ of $n$,
\begin{equation}\label{eq:Oarray-inner-product}
\left\langle
S_O,
 h_{\alpha_1}\cdots h_{\alpha_\ell}
\right\rangle
=
|\Acal_O(\alpha)|.
\end{equation}
\end{lemma}

\begin{proof}
Let $\alpha^+$ be the partition obtained by deleting the zero parts of $\alpha$ and sorting the
remaining parts in weakly decreasing order.  Since multiplication is commutative,
\[
h_{\alpha_1}\cdots h_{\alpha_\ell}=h_{\alpha^+}.
\]
By Hall duality \eqref{eq:Hall-duality},
\[
\langle S_O,h_{\alpha^+}\rangle
=[m_{\alpha^+}]S_O.
\]
Because $S_O$ is symmetric, this is also the coefficient of the particular monomial with exponent
sequence $\alpha$:
\begin{equation}\label{eq:monomial-colored}
[m_{\alpha^+}]S_O
=
[x_1^{\alpha_1}\cdots x_\ell^{\alpha_\ell}]S_O.
\end{equation}
The coefficient in \eqref{eq:monomial-colored} counts pairs $(M,\kappa)$, where $M\in\Mcal(O)$ and
$\kappa$ is a proper coloring of $\Inc(P_M)$ using color $i$ exactly $\alpha_i$ times.

Fix a color $i$.  Two matching edges of the same color are nonadjacent in $\Inc(P_M)$ and hence
comparable in $P_M$.  Thus the color-$i$ edges form a chain, which can be written uniquely as
\[
\{a_1,b_1\}<_{P_M}\cdots<_{P_M}\{a_{\alpha_i},b_{\alpha_i}\}.
\]
By the definition of $P_M$,
\[
a_1<b_1<a_2<b_2<\cdots<a_{\alpha_i}<b_{\alpha_i}.
\]
Since $\Open(M)=O$, every $a_j$ lies in $O$ and every $b_j$ lies outside $O$.  Hence the color class
determines the row
\[
R_i=(a_1,b_1,\ldots,a_{\alpha_i},b_{\alpha_i}),
\]
and collecting the rows over all colors gives an $O$-array of type $\alpha$.

Conversely, take $R\in\Acal_O(\alpha)$ and pair consecutive entries in each row.  Since every
element of $[2n]$ appears exactly once, these pairs form a perfect matching $M$.  Across all rows
there are
\[
\alpha_1+\cdots+\alpha_\ell=n
\]
odd row positions.  Their entries are distinct elements of $O$, and $|O|=n$, so the set of
odd-position entries is exactly $O$.  Hence $\Open(M)=O$.

Color every edge obtained from row $R_i$ by color $i$.  Since each row is strictly increasing,
consecutive paired edges satisfy
\[
\{r_{i,2j-1},r_{i,2j}\}<_{P_M}
\{r_{i,2j+1},r_{i,2j+2}\}.
\]
Thus each color class is a chain in $P_M$, equivalently a stable set in $\Inc(P_M)$, so the coloring
is proper.  The two constructions are inverse to one another, proving
\eqref{eq:Oarray-inner-product}.
\end{proof}

\subsection{Jacobi--Trudi cancellation}

We prove \cref{thm:fixedO} coefficientwise.  Fix a partition
\[
\mu=(\mu_1,\ldots,\mu_\ell)\vdash n.
\]
Transposing the matrix in the ordinary Jacobi--Trudi identity \cite[Theorem~7.16.1]{StanleyEC2} gives
\begin{equation}\label{eq:transposed-JT}
\begin{aligned}
s_\mu
&=
\det\bigl(h_{\mu_j-j+i}\bigr)_{1\le i,j\le\ell}\\
&=
\sum_{\pi\in \mathfrak S_\ell}\sgn(\pi)
\prod_{i=1}^\ell h_{\alpha_i(\pi)},
\end{aligned}
\end{equation}
where
\begin{equation}\label{eq:alpha-pi}
\alpha_i(\pi):=\mu_{\pi(i)}-\pi(i)+i.
\end{equation}
Terms with some $\alpha_i(\pi)<0$ vanish by the convention $h_r=0$ for $r<0$;
we correspondingly set $\Acal_O(\alpha(\pi))=\varnothing$.

Define the signed set
\begin{equation}\label{eq:signed-set}
\frakA_O(\mu)
:=
\bigsqcup_{\pi\in \mathfrak S_\ell}
\{\pi\}\times\Acal_O(\alpha(\pi)),
\end{equation}
where an element $(\pi,R)$ is assigned sign $\sgn(\pi)$.

For $R=(R_1,\ldots,R_\ell)$ and $t\in\{0,1,\ldots,2n\}$, define the prefix count
\begin{equation}\label{eq:prefix-count}
N_i(t):=\#\{x\in R_i:x\le t\}.
\end{equation}
Call $R$ \emph{good} if
\begin{equation}\label{eq:good-condition}
N_{i+1}(t)\le N_i(t)+1
\end{equation}
for every $i\in[\ell-1]$ and every $t\in[2n]$; otherwise call it \emph{bad}.

The involution below is modeled on the array involutions used by Gasharov \cite{Gasharov} and
restated by Chow \cite[Chapter~2, Section~4, Theorem~5]{ChowThesis}.  The parity argument preserving
the prescribed opening set is the feature specific to the present refinement.

\begin{lemma}[Sign-reversing cancellation]\label{lem:involution}
There is a sign-reversing involution on the bad elements of $\frakA_O(\mu)$.  The uncancelled good
elements are naturally in bijection with
\[
\{A\in\SYT(\rho(\mu')):\Left(A)=O\}.
\]
Consequently,
\begin{equation}\label{eq:signed-cancel}
\sum_{\pi\in \mathfrak S_\ell}
\sgn(\pi)|\Acal_O(\alpha(\pi))|
=
g_O^{\mu'}.
\end{equation}
\end{lemma}

\begin{proof}
Take a bad element $(\pi,R)$.  Choose the smallest integer $t$ for which \eqref{eq:good-condition}
fails for some adjacent pair of rows.

Let $R_q$ be the row containing $t$. Then $N_q(t)=N_q(t-1)+1$, while
$N_i(t)=N_i(t-1)$ for every $i\ne q$. Set
\[
\Delta_i(t):=N_{i+1}(t)-N_i(t)\qquad(1\le i<\ell).
\]
Among the defined adjacent differences, $\Delta_{q-1}$ increases by one if $q>1$, while $\Delta_q$
decreases by one if $q<\ell$; all others are unchanged. Since no difference exceeded $1$ before the
first bad value, a new violation requires $q>1$ and can occur only in $\Delta_{q-1}$, which must
then equal $2$. Thus the failing adjacent pair at
the first bad value is unique.  Writing $r=q-1$, we have
\begin{equation}\label{eq:first-failure}
N_{r+1}(t)=N_r(t)+2.
\end{equation}

Set $p:=N_r(t)$.  Write
\[
R_r=(u_1<\cdots<u_p<u_{p+1}<\cdots<u_{2a}),
\qquad
a=\alpha_r(\pi),
\]
and
\[
R_{r+1}=(v_1<\cdots<v_{p+2}<v_{p+3}<\cdots<v_{2b}),
\qquad
b=\alpha_{r+1}(\pi).
\]
Here the first $p$ entries of $R_r$ and the first $p+2$ entries of $R_{r+1}$ are precisely
the entries at most $t$ in those rows.

Exchange the two tails above $t$ by setting
\begin{align}
R'_r
&=(u_1,\ldots,u_p,v_{p+3},\ldots,v_{2b}),
\label{eq:tail-swap-1}\\
R'_{r+1}
&=(v_1,\ldots,v_{p+2},u_{p+1},\ldots,u_{2a}),
\label{eq:tail-swap-2}
\end{align}
and leave all other rows unchanged.  At the same time, let
\begin{equation}\label{eq:pi-prime}
\pi':=\pi s_r,
\qquad
s_r=(r\ r+1).
\end{equation}
Thus $s_r$ is the adjacent transposition, and right multiplication by $s_r$ swaps the values of
$\pi$ in positions $r$ and $r+1$:
\[
\pi'(r)=\pi(r+1),
\qquad
\pi'(r+1)=\pi(r).
\]

Every prefix entry is at most $t$, whereas every transferred tail entry is greater than $t$, so the
new rows remain strictly increasing.  No entry is lost or repeated.  Moreover, the cuts occur after
positions $p$ and $p+2$, which have the same parity.  Therefore every transferred entry keeps the
parity of its row position.  Odd positions remain occupied by elements of $O$, and even positions
remain occupied by elements outside $O$; hence the alternating opener/closer condition is preserved.

The first failure implies $p+2\le2b$, so $b\ge1$. The old half-row-lengths are $(a,b)$.  The new row
lengths are
\[
|R'_r|=2(b-1),
\qquad
|R'_{r+1}|=2(a+1),
\]
so
\begin{equation}\label{eq:half-row-change}
(a,b)\longmapsto(b-1,a+1).
\end{equation}
On the other hand, from \eqref{eq:alpha-pi} and \eqref{eq:pi-prime},
\[
\alpha_r(\pi')
=\alpha_{r+1}(\pi)-1=b-1,
\qquad
\alpha_{r+1}(\pi')
=\alpha_r(\pi)+1=a+1.
\]
Thus $R'\in\Acal_O(\alpha(\pi'))$.

Since an adjacent transposition is odd,
\begin{equation}\label{eq:sign-reversal}
\sgn(\pi')
=
\sgn(\pi)\sgn(s_r)
=
-\sgn(\pi).
\end{equation}
No entry at most $t$ has moved, so for every $s\le t$ and every row $i$,
\[
N'_i(s)=N_i(s).
\]
Consequently the same first bad value $t$ and the same adjacent pair $r,r+1$ are selected after the
first tail exchange.  A second application swaps the same tails back, and
\[
\pi's_r=\pi s_rs_r=\pi.
\]
The map is therefore a sign-reversing involution on the bad elements of $\frakA_O(\mu)$.

It remains to identify the survivors.  Suppose $(\pi,R)$ is good.  At $t=2n$ we have
$N_i(2n)=2\alpha_i(\pi)$, so \eqref{eq:good-condition} gives
\[
2\alpha_{i+1}(\pi)\le2\alpha_i(\pi)+1.
\]
Since both doubled quantities are even,
\begin{equation}\label{eq:alpha-weakly-decreasing}
\alpha_{i+1}(\pi)\le\alpha_i(\pi).
\end{equation}
Let $c_j:=\mu_j-j$.  Because $\mu$ is a partition,
\[
c_1>c_2>\cdots>c_\ell.
\]
Using $\alpha_i(\pi)=c_{\pi(i)}+i$, inequality \eqref{eq:alpha-weakly-decreasing} implies
\[
c_{\pi(i)}>c_{\pi(i+1)},
\]
and hence $\pi(i)<\pi(i+1)$ for every $i$.  Therefore
\begin{equation}\label{eq:pi-id}
\pi=\mathrm{id}.
\end{equation}
The surviving row lengths are thus $2\mu_1,\ldots,2\mu_\ell$.

Write
\[
R_i=(r_{i,1}<\cdots<r_{i,2\mu_i}).
\]
We claim that the prefix condition \eqref{eq:good-condition} is equivalent to
\begin{equation}\label{eq:column-condition}
r_{i,j}<r_{i+1,j+1}
\end{equation}
whenever both entries exist.  If \eqref{eq:column-condition} fails, then with $t=r_{i+1,j+1}$ one gets
\[
N_{i+1}(t)\ge j+1,
\qquad
N_i(t)\le j-1,
\]
contradicting goodness.  Conversely, if goodness fails at $t$, set $k=N_{i+1}(t)$.  Since
$\pi=\mathrm{id}$ and $\mu_i\ge\mu_{i+1}$, the entry $r_{i,k-1}$ exists and
\[
r_{i+1,k}\le t<r_{i,k-1},
\]
which violates \eqref{eq:column-condition} with $j=k-1$.

The row lengths of $\rho(\mu')$ are exactly $2\mu_1,\ldots,2\mu_\ell$, and the $j$th cell of row $i$
is vertically aligned with the $(j+1)$st cell of row $i+1$.  Hence \eqref{eq:column-condition} is
precisely the column-increasing condition for a standard Young tableau of shape $\rho(\mu')$.  The
rows are already strictly increasing.  Finally, the left cell of each horizontal domino is an odd
row position, and the $O$-array condition says that the set of entries in all odd row positions is
exactly $O$.  Thus the survivors are precisely the tableaux in
\[
\{A\in\SYT(\rho(\mu')):\Left(A)=O\}.
\]
This establishes \eqref{eq:signed-cancel}.
\end{proof}

\begin{proof}[Proof of \Cref{thm:fixedO}]
By Schur orthonormality in \eqref{eq:Hall-duality} and \eqref{eq:transposed-JT},
\begin{align*}
[s_\mu]S_O
&=\langle S_O,s_\mu\rangle\\
&=
\sum_{\pi\in \mathfrak S_\ell}
\sgn(\pi)
\left\langle
S_O,
\prod_{i=1}^\ell h_{\alpha_i(\pi)}
\right\rangle.
\end{align*}
By \cref{lem:Oarray},
\[
[s_\mu]S_O
=
\sum_{\pi\in \mathfrak S_\ell}
\sgn(\pi)|\Acal_O(\alpha(\pi))|.
\]
By \cref{lem:involution}, the sum on the right-hand side equals $g_O^{\mu'}$.
Since $\mu\vdash n$ was arbitrary,
\[
S_O
=
\sum_{\mu\vdash n}g_O^{\mu'}s_\mu
=
\sum_{\lambda\vdash n}g_O^\lambda s_{\lambda'}.
\]
Applying $\omega$ and using $\omega(s_{\lambda'})=s_\lambda$ gives \eqref{eq:fixedO-schur-omega}.
\end{proof}

Every perfect matching has a unique opening set, and every tableau $A\in\SYT(\rho(\lambda))$ has a
unique left set $\Left(A)$.  Consequently, summing over all valid opening sets $O\subseteq[2n]$,
\[
\sum_O S_O
=
\sum_{M\in\Mcal(2n)}X_{\Inc(P_M)},
\qquad
\sum_O g_O^\lambda=f^{\rho(\lambda)}.
\]
Here $f^\theta=|\SYT(\theta)|$ for a skew shape $\theta$. Thus \cref{thm:fixedO} refines the global Schur expansion in \cite{ASS} fiber by fiber.

\subsection{Descent-set refinements}\label{sec:cardinality}

Chow's sequencing expansion \cite[Corollary~2]{Chow} gives
\begin{equation}\label{eq:Chow-sequencing}
 X_{\Inc(P_M)}=\sum_{\sigma}Q_{\DES(\sigma),n},
\end{equation}
where $\sigma$ ranges over sequencings of the edges of $M$. Summing over $M\in\Mcal(O)$ yields
\begin{equation}\label{eq:SO-Q}
 S_O=\sum_{D\subseteq[n-1]}N_O(D)Q_{D,n}.
\end{equation}
On the other hand, the standard fundamental expansion is
\begin{equation}\label{eq:Schur-Q}
 s_\mu=\sum_{D\subseteq[n-1]}f_D^\mu Q_{D,n};
\end{equation}
see \cite[Section~7.19]{StanleyEC2}. We also use the following standardization identity.

\begin{lemma}[Standardization]\label{lem:standardization}
For $\mu\vdash n$ and $T\subseteq[n-1]$,
\begin{equation}\label{eq:standardization}
 K_{\mu,\alpha(T)}=\sum_{D\subseteq T}f_D^\mu.
\end{equation}
\end{lemma}
\begin{proof}
Write $T=\{t_1<\cdots<t_k\}$, with $t_0=0$ and $t_{k+1}=n$. In a semistandard tableau of content
$\alpha(T)$, replace the copies of $r$, from left to right, by $t_{r-1}+1,\ldots,t_r$. Equal entries
occupy distinct columns and, in increasing column order, lie weakly higher. The resulting tableau is
standard and has no descent inside any of these blocks. Conversely, in a standard tableau with
descent set contained in $T$, replace the labels $t_{r-1}+1,\ldots,t_r$ by $r$. Within such a block
the row indices weakly decrease, so two labels of the block cannot occupy the same column. The
result is semistandard. Consecutive labels in a block without descents also move strictly to the
right in column order, so standardization recovers the original labels. The constructions are
inverse.
\end{proof}

\begin{corollary}[Descent-set refinements]\label{thm:cardinality}\label{cor:exact-descent}
For every valid opening set $O\subseteq[2n]$ and all $D,T\subseteq[n-1]$,
\begin{align}
 N_O(D)&=\sum_{\lambda\vdash n}g_O^\lambda f_D^{\lambda'},
       \label{eq:exact-descent}\\
 R_O(T)&=\sum_{\lambda\vdash n}g_O^\lambda K_{\lambda',\alpha(T)}.
       \label{eq:cardinality-identity}
\end{align}
\end{corollary}
\begin{proof}
Compare the coefficients of $Q_{D,n}$ in \eqref{eq:SO-Q} and in the expansion obtained from
\cref{thm:fixedO} and \eqref{eq:Schur-Q}. This gives \eqref{eq:exact-descent}. Summing over
$D\subseteq T$ and applying \cref{lem:standardization} gives \eqref{eq:cardinality-identity}.
\end{proof}

In particular, define
\begin{align}
 \Ccal_{O,T}
 &=\{(M,\sigma):M\in\Mcal(O),\ \DES(\sigma)\subseteq T\},\label{eq:COT}\\
 \Tcal_{O,T}
 &=\bigsqcup_{\lambda\vdash n}
 \bigl(\{A\in\SYT(\rho(\lambda)):\Left(A)=O\}
       \times\SSYT(\lambda',\alpha(T))\bigr).\label{eq:TOT}
\end{align}
Then \eqref{eq:cardinality-identity} says that
\begin{equation}\label{eq:set-cardinality}
 |\Ccal_{O,T}|=|\Tcal_{O,T}|.
\end{equation}
The coefficient comparison gives this equality without constructing a bijection between the two sets.

\section{Selected coefficients and sink refinements}\label{sec:hcoeff}

We now give direct formulas for several coefficients in
\begin{equation}\label{eq:HO-h}
H_O=\sum_{\lambda\vdash n}c^O_\lambda h_\lambda.
\end{equation}
For a graph $G$, let $\sink(G,k)$ denote the number of acyclic orientations of $G$ having exactly $k$ sinks.
An orientation assigns a direction to each edge; it is acyclic if it has no directed cycle, and
a sink is a vertex with no outgoing edge.

\subsection{Aggregation by the number of sinks}

\begin{proposition}[Sink identity]\label{thm:fixedOsink}
For every valid opening set $O\subseteq[2n]$ and every $1\le k\le n$,
\begin{equation}\label{eq:fixedO-sink}
\sum_{\substack{\lambda\vdash n\\\ell(\lambda)=k}}
 c^O_\lambda
=
\sum_{M\in\Mcal(O)}
\sink(\Inc(P_M),k).
\end{equation}
\end{proposition}

\begin{proof}
For each $M$, Stanley's sink theorem \cite[Theorem~3.3]{StanleyXG} gives
\[
 \sink(\Inc(P_M),k)
 =\sum_{\substack{\lambda\vdash n\\\ell(\lambda)=k}}
   [e_\lambda]X_{\Inc(P_M)}.
\]
Sum over $M\in\Mcal(O)$ and use $[e_\lambda]S_O=[h_\lambda]H_O=c^O_\lambda$.
\end{proof}

The lengths $n$, $n-1$ (for $n\ge2$), and $1$ each specify a unique partition of $n$, so these cases
determine individual coefficients.

\subsection{The coefficient of \texorpdfstring{$h_1^n$}{h1 to the n}}

Define
\[
O_{\min}:=\{1,3,5,\ldots,2n-1\},
\qquad
M_{\min}:=12\mid34\mid56\mid\cdots\mid(2n-1)(2n).
\]

\begin{theorem}\label{thm:h1n}
For every valid opening set $O\subseteq[2n]$,
\begin{equation}\label{eq:h1n}
[h_1^n]H_O
=
\begin{cases}
1,&O=O_{\min},\\
0,&O\ne O_{\min}.
\end{cases}
\end{equation}
Moreover,
\begin{equation}\label{eq:Homin}
H_{O_{\min}}=h_1^n.
\end{equation}
\end{theorem}

\begin{proof}
The partition $(1^n)$ is the unique partition of $n$ with $n$ parts.  Setting $k=n$ in
\cref{thm:fixedOsink} gives
\begin{equation}\label{eq:h1n-sink}
[h_1^n]H_O
=
\sum_{M\in\Mcal(O)}\sink(\Inc(P_M),n).
\end{equation}
An acyclic orientation of an $n$-vertex graph has $n$ sinks if and only if the graph is edgeless.
Thus the right-hand side of \eqref{eq:h1n-sink} counts those $M\in\Mcal(O)$ for which $\Inc(P_M)$ is
edgeless.

If $\Inc(P_M)$ is edgeless, all matching intervals are pairwise comparable in $P_M$, hence disjoint
and totally ordered from left to right.  Since the endpoints are $1,2,\ldots,2n$, the leftmost
interval must be $[1,2]$.  Removing it and repeating the argument yields inductively
\[
M=M_{\min}.
\]
Hence $\Open(M)=O_{\min}$.  Conversely, the intervals of $M_{\min}$ are pairwise disjoint, so
$\Inc(P_{M_{\min}})$ is edgeless.  This proves \eqref{eq:h1n}.

For $O=O_{\min}$, a left-to-right scan forces the closer $2i$ to be paired with $2i-1$: by
induction, all earlier openers have already been used. Thus $\Mcal(O_{\min})=\{M_{\min}\}$.  For the
edgeless graph on $n$ vertices,
\[
X_{\Inc(P_{M_{\min}})}=(x_1+x_2+\cdots)^n=p_1^n=h_1^n.
\]
Since $\omega(h_1)=h_1$, equation \eqref{eq:Homin} follows.
\end{proof}

\subsection{The coefficient of \texorpdfstring{$h_2h_1^{n-2}$}{h2 times h1 to the (n-2)}}

Recall the opener--closer word $w(O)$ from \cref{subsec:matchings}. For $1\le i<j\le n$, define $O_{i,j}$ by
\begin{equation}\label{eq:Oij-word}
w(O_{i,j})
=(OC)^{i-1}\,OO\,(CO)^{j-i-1}\,CC\,(OC)^{n-j}.
\end{equation}
This is a Dyck word, so $O_{i,j}$ is a valid opening set.

\begin{theorem}\label{thm:h2h1}
Assume $n\ge2$.  For every valid opening set $O\subseteq[2n]$,
\begin{equation}\label{eq:h2h1}
[h_2h_1^{n-2}]H_O
=
\begin{cases}
4,&O=O_{i,j}\text{ for some }1\le i<j\le n,\\
0,&\text{otherwise}.
\end{cases}
\end{equation}
The pair $(i,j)$ is uniquely determined by $O_{i,j}$.  Consequently,
\begin{equation}\label{eq:h2h1-global}
\sum_O[h_2h_1^{n-2}]H_O
=4\binom n2
=2n(n-1).
\end{equation}
\end{theorem}

\begin{proof}
The partition $(2,1^{n-2})$ is the unique partition of $n$ having $n-1$ parts.  By \cref{thm:fixedOsink},
\begin{equation}\label{eq:h2h1-sink}
[h_2h_1^{n-2}]H_O
=
\sum_{M\in\Mcal(O)}\sink(\Inc(P_M),n-1).
\end{equation}

An $n$-vertex graph $G$ admits an acyclic orientation with exactly $n-1$ sinks if and only if it has
at least one edge and has a vertex incident with every edge. More precisely, such orientations are
in bijection with choices of a vertex $v$ incident with every edge.  Indeed, if every vertex other
than $v$ is a sink, no two of those sinks can be adjacent, so every edge is incident with $v$.
Conversely, if every edge is incident with $v$ and at least one edge exists, orient every edge away
from $v$; this is the unique orientation in which all other vertices are sinks.

Therefore the right-hand side of \eqref{eq:h2h1-sink} counts pairs $(M,v)$ such that $M\in\Mcal(O)$
and $\Inc(P_M)$ is a nonempty star centered at $v$, together with isolated vertices.  Regard $v$ as
the matching edge
\[
v=\{a,b\},\qquad a<b.
\]
After removing $v$, all remaining matching intervals are pairwise comparable.  If
\[
[2n]\setminus\{a,b\}=\{d_1<d_2<\cdots<d_{2n-2}\},
\]
then the remaining matching is forced to be
\begin{equation}\label{eq:forced-pairing}
\{d_1,d_2\},\{d_3,d_4\},\ldots,
\{d_{2n-3},d_{2n-2}\}.
\end{equation}
Conversely, for any $a<b$, pairing the remaining points as in \eqref{eq:forced-pairing} ensures that
every incomparability edge is incident with $v$.

Partition $[2n]$ into consecutive blocks
\[
I_t:=\{2t-1,2t\},\qquad1\le t\le n.
\]
Let $a\in I_i$ and $b\in I_j$.  Since $a<b$, we have $i\le j$.  If $i=j$, then $\{a,b\}=I_i$ and
\eqref{eq:forced-pairing} gives $M_{\min}$, whose incomparability graph is edgeless.  If $i<j$, the
distinguished edge is not an edge of $M_{\min}$; by the edgeless characterization in \cref{thm:h1n},
the incomparability graph has at least one edge. Hence the pairs counted by \eqref{eq:h2h1-sink} are exactly
those with
\begin{equation}\label{eq:i-less-j}
i<j.
\end{equation}

Fix $i<j$.  There are two choices for $a\in I_i$ and two choices for $b\in I_j$, hence four choices
for the distinguished edge $v=\{a,b\}$.  The opening set of the resulting matching is determined as
follows.  For $x\in[2n]\setminus\{a,b\}$, the point $x$ is an opener precisely when its rank in the
increasing list $[2n]\setminus\{a,b\}$ is odd.  Equivalently,
\begin{equation}\label{eq:opener-rank}
x\text{ is an opener}
\Longleftrightarrow
\begin{cases}
x\text{ is odd},&x<a,\\
x\text{ is even},&a<x<b,\\
x\text{ is odd},&x>b.
\end{cases}
\end{equation}
Including the opener $a$ itself, the block pattern is exactly
\[
(OC)^{i-1}\,OO\,(CO)^{j-i-1}\,CC\,(OC)^{n-j}.
\]
Thus all four choices give the same opening set $O_{i,j}$.

Conversely, every pair $(M,v)$ counted by \eqref{eq:h2h1-sink} arises from a unique pair of blocks
$I_i,I_j$ with $i<j$ and one of the four choices $a\in I_i$, $b\in I_j$.  In \eqref{eq:Oij-word},
the unique block $OO$ occurs at $i$ and the unique block $CC$ occurs at $j$, so distinct pairs
$(i,j)$ give distinct opening sets.  Therefore the coefficient is $4$ for $O=O_{i,j}$ and $0$
otherwise.  These are four distinct \emph{marked pairs} $(M,v)$, not necessarily four distinct
matchings; for a single-edge incomparability graph the two choices of center correspond to its two
orientations. Summing over the $\binom n2$ pairs $i<j$ proves \eqref{eq:h2h1-global}.
\end{proof}

\subsection{The coefficient of \texorpdfstring{$h_n$}{hn}}

Recall $C=\{c_1<\cdots<c_n\}$ and the pre-closure heights $\eta_i(O)$ from \eqref{eq:eta}. We have
$\eta_i(O)\ge1$ and $\eta_n(O)=1$.

Let $\chi_G(q)$ denote the chromatic polynomial of $G$. For an integer $q\ge0$, it counts proper
colorings with colors in $[q]$. We write $[q]\chi_G(q)$ for its linear coefficient.

\begin{lemma}\label{lem:fixedOchromatic}
For every $M\in\Mcal(O)$,
\begin{equation}\label{eq:fixedOchromatic}
\chi_{\Inc(P_M)}(q)
=
\prod_{i=1}^n\bigl(q-(\eta_i(O)-1)\bigr).
\end{equation}
In particular, the chromatic polynomial of $\Inc(P_M)$ depends only on $O$, not on the choice of
$M\in\Mcal(O)$.
\end{lemma}

\begin{proof}
For each $i$, let $v_i$ be the matching interval whose right endpoint is $c_i$.  Order the vertices as
\[
v_1,v_2,\ldots,v_n.
\]
Immediately before $c_i$ closes, exactly $\eta_i(O)$ matching intervals are open.  One is $v_i$, and
the remaining $\eta_i(O)-1$ close later.  These are precisely the later neighbors of $v_i$ among
$v_{i+1},\ldots,v_n$.

All these later intervals contain the point $c_i$, so they are pairwise incomparable in $P_M$, hence
pairwise adjacent in $\Inc(P_M)$.  Therefore the later neighbors of $v_i$ form a clique of size
$\eta_i(O)-1$.

Fix an integer $q\ge n$ and color in the order $v_n,v_{n-1},\ldots,v_1$.  When $v_i$ is colored, its
already colored neighbors form a clique and therefore use $\eta_i(O)-1$ distinct colors.  Hence
there are
\[
q-(\eta_i(O)-1)
\]
choices for the color of $v_i$.  Multiplying over $i$ proves \eqref{eq:fixedOchromatic} for all
integers $q\ge n$. Both sides are polynomials in $q$, so the identity holds as a polynomial
identity.
\end{proof}

\begin{theorem}\label{thm:hn}
For every valid opening set $O\subseteq[2n]$,
\begin{equation}\label{eq:hn}
[h_n]H_O
=
n\prod_{i=1}^{n-1}\eta_i(O)\bigl(\eta_i(O)-1\bigr).
\end{equation}
Consequently,
\begin{equation}\label{eq:hn-positive}
[h_n]H_O>0
\quad\Longleftrightarrow\quad
\eta_i(O)\ge2\text{ for every }i<n.
\end{equation}
Equivalently, $[h_n]H_O>0$ if and only if the Dyck word $w(O)$ does not return to height $0$ before
its final letter.
\end{theorem}

\begin{proof}
The partition $(n)$ is the unique partition of $n$ having one part.  Taking $k=1$ in
\cref{thm:fixedOsink} gives
\begin{equation}\label{eq:hn-sink}
[h_n]H_O
=
\sum_{M\in\Mcal(O)}\sink(\Inc(P_M),1).
\end{equation}
Fix $M\in\Mcal(O)$ and set $G=\Inc(P_M)$.  By the Greene--Zaslavsky unique-sink theorem
\cite[Theorem~7.3]{GreeneZaslavsky}, in the form recorded by Stanley
\cite[Theorem~1.2(b)]{StanleyXG}, for any prescribed vertex the number of acyclic orientations of
$G$ having that vertex as unique sink is
\[
(-1)^{n-1}[q]\chi_G(q).
\]
Since there are $n$ possible choices of the unique sink,
\begin{equation}\label{eq:unique-sink-total}
\sink(G,1)
=
n(-1)^{n-1}[q]\chi_G(q).
\end{equation}
By \cref{lem:fixedOchromatic} and \eqref{eq:eta-basic},
\[
\chi_G(q)
=
q\prod_{i=1}^{n-1}\bigl(q-(\eta_i(O)-1)\bigr).
\]
Thus
\[
(-1)^{n-1}[q]\chi_G(q)
=
\prod_{i=1}^{n-1}\bigl(\eta_i(O)-1\bigr).
\]
Every matching $M\in\Mcal(O)$ therefore contributes the same quantity
\[
n\prod_{i=1}^{n-1}\bigl(\eta_i(O)-1\bigr)
\]
to \eqref{eq:hn-sink}.  Using \cref{lem:matchingheight},
\begin{align*}
[h_n]H_O
&=
|\Mcal(O)|\,
 n\prod_{i=1}^{n-1}\bigl(\eta_i(O)-1\bigr)\\
&=
n\prod_{i=1}^{n-1}\eta_i(O)\bigl(\eta_i(O)-1\bigr),
\end{align*}
which proves \eqref{eq:hn}.  Since all $\eta_i(O)$ are positive integers, \eqref{eq:hn-positive}
follows immediately.

Finally, $\eta_i(O)=1$ for some $i<n$ precisely when the Dyck path returns to height $0$ at the
closer $c_i$.  This gives the last equivalence.
\end{proof}

For the running example, $O=\{1,2,4,6\}=O_{1,4}$ and $\eta(O)=(2,2,2,1)$.
\Cref{thm:h1n,thm:h2h1,thm:hn} therefore give
$[h_1^4]H_O=0$, $[h_2h_1^2]H_O=4$, and $[h_4]H_O=4(2\cdot1)^3=32$.

The one-part coefficient also refines the cyclically alternating interpretation of $[h_n](2n)!A_n$
in \cite[Theorem~5.2(c)]{ASS}.

\begin{definition}\label{def:cyclic-alt}
Let $\CAcal_{2n}(O)$ be the set of permutations
\[
w=a_1a_2\cdots a_{2n}\in \mathfrak S_{2n}
\]
satisfying
\begin{equation}\label{eq:cyclic-alt}
a_1>a_2<a_3>a_4<\cdots>a_{2n-1}>a_{2n},
\qquad
a_{2n}<a_1,
\end{equation}
and whose set of valley values is exactly
\begin{equation}\label{eq:valley-set}
\{a_2,a_4,\ldots,a_{2n}\}=O.
\end{equation}
These are linear permutations; cyclic rotations are not identified.
\end{definition}

\begin{theorem}[Fixed-opening-set cyclic interpretation]\label{thm:hncyclic}
For every valid opening set $O\subseteq[2n]$,
\begin{equation}\label{eq:hncyclic}
[h_n]H_O
=
|\CAcal_{2n}(O)|.
\end{equation}
Thus the coefficient of $h_n$ in $H_O$ counts cyclically alternating permutations whose valley-value
set is exactly $O$.
\end{theorem}

\begin{proof}
Let
\[
C=[2n]\setminus O=\{c_1<c_2<\cdots<c_n\}.
\]
In a permutation satisfying \eqref{eq:cyclic-alt}, the elements of $O$ are valley values and the
elements of $C$ are peak values.  Since $2n$ cannot be a valley, $c_n=2n$.

For $n=1$, both sides of \eqref{eq:hncyclic} equal $1$, so assume $n\ge2$. First count cyclic
alternating words up to cyclic rotation, retaining their direction; reversal is not identified.
Process the labels $1,2,\ldots,2n$ in increasing order while maintaining a collection of oriented
paths whose initial and terminal vertices are valleys.  If the current label $x$ lies in $O$, create
the singleton oriented path $(x)$.  If the current label is $c_i\in C$ with $i<n$, then immediately
before inserting $c_i$ the number of oriented paths is
\[
\#(O\cap[c_i-1])-(i-1)=\eta_i(O).
\]
Choose an ordered pair of distinct paths $(P,Q)$ and concatenate them as
\[
P,\ c_i,\ Q.
\]
There are therefore $\eta_i(O)\bigl(\eta_i(O)-1\bigr)$ choices for the insertion of the peak $c_i$.

At the final peak $c_n=2n$, equation \eqref{eq:eta-basic} gives $\eta_n(O)=1$, so exactly one
oriented path remains.  Joining its terminal valley to its initial valley through the peak $2n$
produces an oriented cyclic alternating arrangement with valley-value set $O$.

To see both exhaustiveness and uniqueness, take any such cyclic word. Its restriction to labels at
most $t<2n$ is a disjoint union of directed paths, including singletons. A valley is introduced as a
singleton because both of its neighbors are larger. A nonfinal peak joins two distinct components:
joining a component to itself would create a cycle before the full word was assembled. Thus an early
height $\eta_i(O)=1$ makes the construction impossible and gives zero on both sides of the counting
formula. Otherwise the construction proceeds as described.

It is reversible. Starting from a cyclic word, first remove its final peak $2n$ to obtain a path,
and then remove the remaining labels in decreasing order.  Removing a peak splits one oriented path
into the ordered pair of paths that had been joined when that peak was inserted.  When a valley is
removed, its two larger neighboring peaks have already been removed, so it appears as a singleton
path.  Therefore the number of oriented cyclic arrangements is
\begin{equation}\label{eq:oriented-cyclic-count}
\prod_{i=1}^{n-1}\eta_i(O)\bigl(\eta_i(O)-1\bigr).
\end{equation}
Each oriented cyclic arrangement yields a linear permutation satisfying \eqref{eq:cyclic-alt} after
choosing which of its $n$ peaks appears in the first position.  Hence
\[
|\CAcal_{2n}(O)|
=
n\prod_{i=1}^{n-1}\eta_i(O)\bigl(\eta_i(O)-1\bigr).
\]
Comparison with \cref{thm:hn} proves \eqref{eq:hncyclic}.
\end{proof}

Let $E_m$ denote the Euler numbers, defined by
\[
 \sec t+\tan t=\sum_{m\ge0}E_m\frac{t^m}{m!}.
\]
Equivalently, $E_m$ counts alternating permutations of $[m]$ with the specified initial direction
$>$; reversing the initial direction gives the same number. These conventions agree with
\cite[Section~3]{ASS}.

\begin{corollary}\label{cor:globalhn}
Summing over all valid opening sets gives
\begin{equation}\label{eq:global-hn}
\sum_O[h_n]H_O
=nE_{2n-1}.
\end{equation}
\end{corollary}

\begin{proof}
The sets $\CAcal_{2n}(O)$ partition the cyclically alternating permutations according to their
valley-value sets.  By \cref{thm:hncyclic},
\[
\sum_O[h_n]H_O
=
\sum_O|\CAcal_{2n}(O)|.
\]
The right-hand side is $nE_{2n-1}$ by \cite[Theorem~5.2(c)]{ASS}.
\end{proof}

\section{Inverse Kostka expansion and one-part rim-hook cancellation}\label{sec:rimhook}

The Schur expansion and the inverse Kostka matrix give a signed tableau model for $[h_\mu]H_O$.

A \emph{rim hook} is a skew Young diagram connected through shared sides and containing no $2\times2$ block.  We use the standard
height convention: the height is one less than the number of rows it occupies.  A rim hook in the
Young diagram of $\lambda$ is \emph{special} if it meets the first column.
Let $\SRHT(\lambda;\mu)$ be the set of \emph{special rim-hook tabloids}: decompositions of $\lambda$
into unlabeled special rim hooks whose multiset of sizes is the multiset of parts of $\mu$.  For $R\in\SRHT(\lambda;\mu)$ put
\[
 \rht(R)=\sum_{B\in R}\bigl(\#\{\text{rows occupied by }B\}-1\bigr),
 \qquad \varepsilon(R)=(-1)^{\rht(R)}.
\]
Our Kostka convention is
\begin{equation}\label{eq:rim-kostka}
 h_\mu=\sum_{\lambda\vdash n}K_{\lambda\mu}s_\lambda.
\end{equation}
With this indexing, the inverse Kostka theorem of E{\u{g}}ecio{\u{g}}lu and Remmel
\cite[Theorem~1]{EgeciogluRemmel} states
\begin{equation}\label{eq:rim-inverse-kostka}
 (K^{-1})_{\mu\lambda}=\sum_{R\in\SRHT(\lambda;\mu)}\varepsilon(R).
\end{equation}
Equal parts of $\mu$ do not produce an additional factorial: the rim hooks are not labeled.

Put
\[
 \mathcal Y_O(\lambda)=\{A\in\SYT(\rho(\lambda)):\Left(A)=O\},
 \qquad
 \mathfrak D_O(\mu)=\bigsqcup_{\lambda\vdash n}\mathcal Y_O(\lambda)\times\SRHT(\lambda;\mu),
\]
and give $(A,R)$ the sign $\varepsilon(R)$. For the running tableau
$A\in\mathcal Y_O((2,2))$ from \cref{ex:running-tableau}, take on the ordinary Young diagram of
$(2,2)$ the rim hook with cells $\{(1,2),(2,1),(2,2)\}$ and the singleton $\{(1,1)\}$. These form a
special rim-hook tabloid $R$ of type $(3,1)$. Thus
$(A,R)\in\mathfrak D_O((3,1))$ has $\rht(R)=1$ and sign $-1$.

\begin{proposition}[Signed rim-hook formula]\label{prop:rim-signed}
For every valid $O\subseteq[2n]$ and every $\mu\vdash n$,
\begin{equation}\label{eq:rim-signed}
 [h_\mu]H_O
 =\sum_{(A,R)\in\mathfrak D_O(\mu)}\varepsilon(R).
\end{equation}
\end{proposition}

\begin{proof}
By \cref{thm:fixedO}, $H_O=\sum_{\lambda}g_O^\lambda s_\lambda$.  Inverting \eqref{eq:rim-kostka} gives
$s_\lambda=\sum_\mu(K^{-1})_{\mu\lambda}h_\mu$, and hence
\[
 [h_\mu]H_O=\sum_{\lambda\vdash n}g_O^\lambda(K^{-1})_{\mu\lambda}.
\]
Now use $g_O^\lambda=|\mathcal Y_O(\lambda)|$ and \eqref{eq:rim-inverse-kostka}.
\end{proof}

The tableau $A$ and rim-hook tabloid $R$ in a pair $(A,R)\in\mathfrak D_O(\mu)$ can be represented
on the same $n$ distinguished dominoes. If
$[\lambda]=\{(i,j):1\le j\le\lambda_i\}$, then row $j$ of $\rho(\lambda)$ contains $\lambda'_j$ distinguished
horizontal dominoes, numbered from left to right as $D_{j,1},\ldots,D_{j,\lambda'_j}$.

The map $(i,j)\mapsto D_{j,i}$ identifies the cells of $\lambda$ with these dominoes: column $j$ of
$\lambda$ and row $j$ of $\rho(\lambda)$ each contain $\lambda'_j$ objects.

Thus $(A,R)\in\mathfrak D_O(\mu)$ can be viewed as a fixed-opening staggered tableau with its
dominoes grouped by a special rim-hook tabloid of type $\mu$. For $\mu=(n)$, the cancellation below
moves whole distinguished dominoes and cancels the negative terms in \eqref{eq:rim-signed}.

For $1\le q\le n$ set
\[
 \lambda^{(q)}=(n-q+1,1^{q-1}).
\]
\begin{corollary}[Hook reduction]\label{cor:rim-hook-reduction}
For every valid opening set $O$,
\begin{equation}\label{eq:rim-hook-reduction}
 [h_n]H_O=\sum_{q=1}^n(-1)^{q-1}g_O^{\lambda^{(q)}}.
\end{equation}
\end{corollary}
\begin{proof}
A type-$(n)$ tabloid consists of one rim hook covering the whole Ferrers diagram. Such a diagram is
a rim hook exactly when it is a hook $\lambda^{(q)}$; the unique rim hook occupies $q$ rows and has
sign $(-1)^{q-1}$. Apply \cref{prop:rim-signed}.
\end{proof}

For $A\in\mathcal Y_O(\lambda^{(q)})$, the top row of $\rho(\lambda^{(q)})$ contains $q$ dominoes.
Number them
$D_1,\ldots,D_q$ from left to right, and number the one-domino rows below it, from top to bottom, by
$D_n,D_{n-1},\ldots,D_{q+1}$.  Write $D_i=(u_i,v_i)$ with $u_i$ in the left cell.  Then
\begin{equation}\label{eq:rim-master}
 \{u_1,\ldots,u_n\}=O,\qquad \{v_1,\ldots,v_n\}=[2n]\setminus O,\qquad u_i<v_i.
\end{equation}
Call $D_1\cdots D_n$ the \emph{master word}, with the \emph{cut} after $D_q$.
For $i<n$, write $\epsilon_i=T$ when $v_i<u_{i+1}$ and $\epsilon_i=V$ when $u_{i+1}<v_i$.

\begin{lemma}[Cut criterion]\label{lem:rim-cut}
The master word $D_1\cdots D_n$ gives a tableau of shape $\rho(\lambda^{(q)})$ if and only if
\begin{equation}\label{eq:rim-cut}
 \epsilon_i=T\ (i<q),\qquad \epsilon_i=V\ (i>q),\qquad u_1<v_n\ \text{ if }q<n.
\end{equation}
When $q<n$, there is no condition on $\epsilon_q$.
\end{lemma}

\begin{proof}
The top-row inequalities are $v_i<u_{i+1}$ for $i<q$.  Since successive one-domino rows are
staggered one cell to
the left, the column inequalities in the tail are $u_{i+1}<v_i$ for $q<i<n$.  If $q<n$, the right
cell of the first
tail domino $D_n$ lies below the leftmost top-row cell, giving $u_1<v_n$. Together with $u_i<v_i$,
these are all the row and column inequalities. For $q<n$, no row or column inequality compares $v_q$ with $u_{q+1}$, so no
condition is imposed on $\epsilon_q$.
\end{proof}

\begin{theorem}[One-part rim-hook involution]\label{thm:rim-onepart}
There is a sign-reversing involution on $\mathfrak D_O((n))$.  Its fixed points all have positive
sign and are naturally
in bijection with the cyclically alternating permutations $\CAcal_{2n}(O)$ of \cref{thm:hncyclic}.
Consequently,
\begin{equation}\label{eq:rim-cyclic}
 \sum_{q=1}^n(-1)^{q-1}g_O^{\lambda^{(q)}}=|\CAcal_{2n}(O)|.
\end{equation}
\end{theorem}

\begin{proof}
For $n=1$, fix the unique object. Assume $n\ge2$.
The rim-hook tabloid is unique once $q$ is fixed, so we act only on the tableau and the cut. If $q<n$ and
$\epsilon_q=T$, move $D_{q+1}$ from the bottom of the one-domino tail to the right end of the top
row, changing the
cut from $q$ to $q+1$.  If $1<q<n$ and $\epsilon_q=V$, move the last top-row domino $D_q$ to a new
bottom tail row,
changing $q$ to $q-1$.  If $q=n$, perform the same latter move with $D_n$.  If $q=1$ and
$\epsilon_1=V$, fix the
object.

The cut criterion proves well-definedness.  In the first case the new top condition at $q$ is exactly
$\epsilon_q=T$; all remaining required inequalities are inherited.  In the second case the new tail
condition is
exactly $\epsilon_q=V$.  When $q=n$, the original top row gives $u_1<v_n$, which is the new boundary
condition.
Every move transfers an entire distinguished domino without reversing its entries, so $\Left(A)=O$
is preserved, and
the rim-hook type remains $(n)$.

The map is involutive.  A move $q\mapsto q+1$ is reversed because the old cut criterion gives
$\epsilon_{q+1}=V$ when $q+1<n$, while the $q=n$ rule reverses the boundary case.  Conversely a move
$q\mapsto q-1$ is reversed because the old criterion gives $\epsilon_{q-1}=T$.  Since the sign at cut $q$ is
$(-1)^{q-1}$, every nonfixed move reverses sign.

The fixed points are exactly the objects with $q=1$ and
\begin{equation}\label{eq:rim-fixed}
 u_{i+1}<v_i\quad(1\le i<n),\qquad u_1<v_n.
\end{equation}
Define
\[
 \Psi(A)=v_n\,u_1\,v_1\,u_2\,v_2\cdots u_{n-1}\,v_{n-1}\,u_n.
\]
By \eqref{eq:rim-master} and \eqref{eq:rim-fixed}, this word is cyclically alternating and its
valley-value set is $O$.
Conversely, from a cyclically alternating word $a_1>a_2<a_3>\cdots>a_{2n}$ with $a_{2n}<a_1$ and
valley set $O$,
set $u_i=a_{2i}$ for $1\le i\le n$, $v_i=a_{2i+1}$ for $1\le i<n$, and $v_n=a_1$.
These inequalities recover a unique fixed tableau.
Thus $\Psi$ is a bijection onto $\CAcal_{2n}(O)$, proving \eqref{eq:rim-cyclic}.
\end{proof}

\section{Further questions}\label{sec:further}

The coefficients $c^O_\mu=[h_\mu]H_O$ are nonnegative integers, but a direct positive combinatorial
interpretation for arbitrary $\mu$ remains to be found. One approach is to construct a
sign-reversing involution on $\mathfrak D_O(\mu)$ with positive fixed points, preserving both the
rim-hook type and the opening set. The one-part case is settled by \cref{thm:rim-onepart}. Sagan and
Lee's rooted special-rim-hook involution \cite{SaganLee} provides related local moves; adapting them
here would also require preservation of standardness in the staggered tableau.

The integral extraction \eqref{eq:fos-integral-coefficient} offers an algebraic formulation of the
same problem. In the transfer formula, $D_O$ comes from permuting openers with equal thresholds,
whereas $\gamma!$ comes from labeling the vertices in the replacement cliques. The appendix encodes
the augmented-monomial coefficients of $S_O$ on a common Ferrers board. A refinement of that model
giving the individual coefficients $c^O_\mu$ would connect these two descriptions.

A second question is to construct a direct bijection between $\Ccal_{O,T}$ and $\Tcal_{O,T}$, and
more finely to realize the exact-descent identity \eqref{eq:exact-descent} bijectively. It would
also be interesting to find a symmetric $q$-refinement compatible with fixed opening sets and the
coefficient transfer.

\appendix
\crefalias{section}{appendix}
\crefname{appendix}{appendix}{appendices}
\Crefname{appendix}{Appendix}{Appendices}
\section{Stable partitions on a common Ferrers board}\label{sec:rookappendix}

We give a fixed-opening-set form of Stanley's stable-partition expansion
\cite[Proposition~2.4]{StanleyXG} and Chow's chain-partition/path-cover correspondence
\cite[Chapter~1, Proposition~2]{ChowThesis}. The stable partitions for all $M\in\Mcal(O)$ are
represented on one Ferrers board; the matching permutation determines a column relabeling and hence
the component types. We retain the threshold data $r,D_O$ from \cref{subsec:thresholds} and
$\Wcal(r),G_a,F_r$ from \cref{prop:regrouping}.

\subsection{Stable partitions and permanents}

Let $\Stab(G)$ be the set of set partitions of $V(G)$ into nonempty stable blocks, and let $\type(\Pi)$ denote
the block-size partition of a set partition $\Pi$. For $v\in[n]$, put
\[
 \pred_\Pi(v)=\max\{u<v:u,v\text{ belong to the same block of }\Pi\},
\]
with $\pred_\Pi(v)=0$ if the set is empty. By \eqref{eq:Ga},
\begin{equation}\label{eq:stable-predecessor}
 \Pi\in\Stab(G_a)
 \quad\Longleftrightarrow\quad
 \pred_\Pi(v)<a_v\quad(v\in[n]).
\end{equation}
Indeed, the earlier vertices in the block of $v$ are all nonadjacent to $v$ exactly when the largest
of them is smaller than $a_v$.

Define
\begin{equation}\label{eq:Nr-partition}
 w_r(\Pi)=\#\{a\in\Wcal(r):\pred_\Pi(v)<a_v\le v\text{ for all }v\}.
\end{equation}
Stanley's expansion then gives
\begin{equation}\label{eq:Fr-stable}
 F_r=\sum_{\Pi\text{ a set partition of }[n]}
       w_r(\Pi)\widetilde m_{\type(\Pi)}.
\end{equation}
Here $\widetilde m_\lambda=(\prod_jm_j(\lambda)!)m_\lambda$, as in \cref{subsec:symmetric-functions}.

For an $n\times n$ matrix $B$, its permanent is
$\perm B=\sum_{\pi\in\mathfrak S_n}\prod_{i=1}^n B_{i,\pi(i)}$.

\begin{proposition}[Permanent formulation]\label{prop:permanent}
For the matrix
\[
 B_r(\Pi)_{i,v}=\ind[\pred_\Pi(v)<r_i\le v],
\]
we have
\begin{equation}\label{eq:permanent-multiplicity}
 \perm B_r(\Pi)=D_Ow_r(\Pi),
\end{equation}
and consequently
\begin{equation}\label{eq:SO-permanent}
 S_O=\sum_{\Pi\text{ a set partition of }[n]}
       \perm B_r(\Pi)\,\widetilde m_{\type(\Pi)}.
\end{equation}
\end{proposition}
\begin{proof}
A permutation $\pi$ contributes to $\perm B_r(\Pi)$ exactly when
\[
 \pred_\Pi(\pi(i))<r_i\le\pi(i)\qquad(i\in[n]).
\]
Equivalently, $\pi\in P_r$ and $a^\pi$ is counted by $w_r(\Pi)$. Each such word has $D_O$ preimages
by \cref{prop:regrouping}, giving \eqref{eq:permanent-multiplicity}. Now use \eqref{eq:Fr-stable}
and $S_O=D_OF_r$.
\end{proof}

\subsection{Rook placements and component types}

Define
\begin{equation}\label{eq:Ferrers-board}
 \Gamma_r=\{(u,i)\in[n]^2:u<r_i\},
\end{equation}
where $u$ indexes rows and $i$ indexes columns. Column $i$ has height $r_i-1$. These heights are
weakly increasing and satisfy $r_i-1<i$, so $\Gamma_r$ is a Ferrers board contained in the staircase
$u<i$. Let $\Rcal(\Gamma_r)$ be its partial nonattacking rook placements, namely subsets of
cells with at most one chosen cell in each row and each column; the empty placement is included.

For $\pi\in P_r$, the full placement
\[
 Q_\pi=\{(\pi(i),i):i\in[n]\}
\]
lies in $[n]^2\setminus\Gamma_r$. Relabeling columns by $i\mapsto\pi(i)$ sends $\Gamma_r$ bijectively to
\begin{equation}\label{eq:relabelled-board}
 \{(u,v)\in[n]^2:u<a_v^\pi\},
\end{equation}
the board of strict comparabilities of $P_{M_\pi}$ under the closer-rank labeling. Thus the
following proposition is the chain-partition/path-cover correspondence on a common board.

\begin{proposition}[Common-board correspondence]\label{thm:rook-bijection}
For every $\pi\in P_r$, the map
\begin{equation}\label{eq:partition-to-rooks}
 \Pi\longmapsto R_{\pi,\Pi}
 =\{(\pred_\Pi(v),\pi^{-1}(v)):\pred_\Pi(v)>0\}
\end{equation}
is a bijection from $\Stab(G_{a^\pi})$ to $\Rcal(\Gamma_r)$.
\end{proposition}
\begin{proof}
If $\Pi$ is stable, \eqref{eq:stable-predecessor} gives $\pred_\Pi(v)<r_{\pi^{-1}(v)}$, so every
specified rook lies on $\Gamma_r$. Column coordinates are distinct. Within a block $b_1<\cdots<b_k$,
the positive predecessors are $b_1,\ldots,b_{k-1}$, each used once; distinct blocks are disjoint.
Thus row coordinates are distinct too.

Conversely, replace every rook $(u,i)\in R$ by the arrow
\begin{equation}\label{eq:rook-arrow}
 u\longrightarrow\pi(i)
\end{equation}
on $[n]$. Since $u<r_i\le\pi(i)$, every arrow is strictly increasing. Nonattacking rooks give
indegree and outdegree at most one, so the components are directed paths and isolated vertices. Let
$\Pi_R$ be their vertex-set partition. Along a path $b_1\to\cdots\to b_k$, we have
\[
 \pred_{\Pi_R}(b_j)=b_{j-1}<r_{\pi^{-1}(b_j)}=a^\pi_{b_j}
 \qquad(2\le j\le k).
\]
The predecessor of a first vertex is zero, so \eqref{eq:stable-predecessor} shows that $\Pi_R$ is
stable. Each arrow records consecutive elements of one block. Hence both constructions recover their
starting object.
\end{proof}

Let $\lambda(\pi,R)$ be the decreasing rearrangement of the vertex counts of the components in
\eqref{eq:rook-arrow}, including isolated vertices as parts of size one. The standard
stable-partition expansion and \cref{thm:rook-bijection} give the following formula.

\begin{corollary}[Rook expansion]\label{cor:twisted-rooks}
For every valid opening set $O$,
\begin{equation}\label{eq:twisted-rooks}
 S_O=\sum_{\pi\in P_r}\ \sum_{R\in\Rcal(\Gamma_r)}
       \widetilde m_{\lambda(\pi,R)}.
\end{equation}
\end{corollary}
\begin{proof}
For fixed $\pi$, the correspondence preserves block sizes, so
$X_{G_{a^\pi}}=\sum_{R\in\Rcal(\Gamma_r)}\widetilde m_{\lambda(\pi,R)}$.
Sum over $\pi\in P_r$ and apply \cref{prop:regrouping}.
\end{proof}

There is also a bipartite description. Regard $Q_\pi\cup R$ as a graph on two disjoint copies of
$[n]$, with a cell $(u,i)$ joining row vertex $u$ to column vertex $i$.
Every vertex lies on one $Q_\pi$-edge and at most one $R$-edge. The components are alternating paths: contracting the $Q_\pi$-edges gives
\eqref{eq:rook-arrow}, which has no directed cycles. A part $k$ of $\lambda(\pi,R)$ therefore counts
$k$ edges of $Q_\pi$, or equivalently $2k$ vertices in the corresponding bipartite component.

\begin{example}\label{ex:rook-types}
Return to $n=4$ and $O=\{1,2,4,6\}$ from \cref{ex:running-matching}. Then
\[
 r=(1,1,2,3),\qquad
 \Gamma_r=\{(1,3),(1,4),(2,4)\}.
\]
Take $R=\{(1,3),(2,4)\}$. Both permutations $\pi_0=1234$ and $\pi_1=1324$, written in one-line
notation, lie in $P_r$. For $\pi_0$, the arrows are $1\to3$ and $2\to4$, so
$\lambda(\pi_0,R)=(2,2)$. For $\pi_1$, they are $1\to2$ and $2\to4$, with vertex $3$ isolated, so
$\lambda(\pi_1,R)=(3,1)$. Thus the rook placement alone does not determine its component type.
\end{example}

To compare with Kim, Lee, and Yoo \cite[Sections~2.4 and~3.1]{KimLeeYoo}, read $(u,i)$ as column
$u$, row $i$ in their convention. Since $u<r_i$ if and only if $B_u<i$, these cells lie above the
Dyck path with column heights $(B_1,\ldots,B_n)$. For $\pi=\mathrm{id}$, the arrows are $u\to i$,
and their path decomposition is precisely the linked-rook decomposition; the corresponding extended
placement has type $\lambda(\mathrm{id},R)$, including the isolated vertices. For general $\pi$,
rooks $(u,i)$ and $(u',i')$ are linked when $\pi(i)=u'$, after the relabeling
\eqref{eq:relabelled-board}. Kim, Lee, and Yoo use linked rooks to obtain a $q$-weighted
Hall--Littlewood expansion for a fixed natural unit interval order \cite[Theorem~3.7]{KimLeeYoo}.
Formula \eqref{eq:twisted-rooks} keeps the augmented-monomial basis and lets the column relabeling
vary with the matching.

\section*{Acknowledgments}
The author is grateful to John Shareshian for helpful discussions, guidance, and encouragement.

\section*{Use of AI}
ChatGPT assisted with checking selected proof steps and algebraic and matrix calculations, as well
as with \LaTeX{} preparation, bibliography formatting, and revision of the exposition. The author
takes full responsibility for the content of the manuscript.

\bibliographystyle{amsalpha}
\bibliography{fixed_opening_set_positivity_arxiv}

\end{document}